\documentclass{article}

\usepackage[margin=1in]{geometry}
\usepackage{amsmath,graphicx,float,amsfonts,amssymb,xcolor,bm,booktabs,multirow,caption,amsthm,authblk}
\newcommand{\bR}{\mathbb{R}}
\newcommand{\bP}{\mathbb{P}}
\newcommand{\bE}{\mathbb{E}}
\newcommand{\bZ}{\mathbb{Z}}
\newcommand{\pSamp}{\mathbb{P}^N}
\newcommand{\pxi}{\mathbb{P}}

\newcommand{\exclude}[1]{}
\newcommand{\extr}{\mathrm{EXT}}
\newcommand{\cset}{I}
\newcommand{\sset}{J}
\newcommand{\ssets}{\mathcal{J}}
\newcommand{\randw}{\overline{W}}
\newcommand{\randt}{\overline{T}}
\newcommand{\randh}{\overline{h}}

\newtheorem{thm}{Theorem}
\newtheorem{lem}[thm]{Lemma}
\newtheorem{prop}[thm]{Proposition}
\newtheorem{cor}[thm]{Corollary}

\newtheorem{defn}{Definition}

\newtheorem{exmp}{Example}
\newtheorem{assumption}{Assumption}

\begin{document}

\title{On sample average approximation for two-stage stochastic programs without relatively complete recourse
\thanks{The authors dedicate this paper to Shabbir Ahmed. This work was supported by NSF Award CMMI-1634597.}
}


\author{Rui Chen \thanks{rchen234@wisc.edu}\qquad James Luedtke \thanks{jim.luedtke@wisc.edu}
}

\date{\small{Department of Industrial and Systems Engineering, University of Wisconsin-Madison}
}

\maketitle

\begin{abstract}
{\sloppypar
We investigate sample average approximation (SAA) for two-stage stochastic programs without relatively complete recourse, i.e., for problems in which there are first-stage feasible solutions that are not guaranteed to have a feasible recourse action. As a feasibility measure of the SAA solution, we consider the ``recourse likelihood'', which is the probability that the solution has a feasible recourse action. For $\epsilon \in (0,1)$, we demonstrate that the probability that a SAA solution has recourse likelihood below $1-\epsilon$ converges to zero exponentially fast with the sample size. Next, we analyze the rate of convergence of optimal solutions of the SAA to optimal solutions of the true problem for problems with a finite feasible region, such as bounded integer programming problems. For problems with non-finite feasible region, we propose modified ``padded'' SAA problems and demonstrate in two cases that such problems can yield, with high confidence, solutions that are certain to have a feasible recourse decision. Finally, we conduct a numerical study on a two-stage resource planning problem that illustrates the results, and also suggests there may be room for improvement in some of the theoretical analysis.
}
\end{abstract}

\section{Introduction}
\label{intro}
We consider a stochastic optimization problem of the form
\begin{equation}\label{real}
\min_{x\in X} f(x):=\bE F(x,\bm{\xi}).
\end{equation}
Here $X\subseteq \bR^{n_1-p}\times\bZ^{p}$ is a closed set that represents the feasible region, which can be a mixed-integer set with $p> 0$, $\bm{\xi}$ is a random vector defined on the complete probability space $(\Omega,\Sigma,\pxi)$ with support $\Xi\subseteq\bR^d$ and
$F:\bR^{n_1}\times\bR^d\rightarrow\bR\cup\{+\infty\}$ is an extended-real-valued function. We assume that $F(x,\cdot)$
is measurable for each $x \in X$.  In particular, we focus on two-stage
stochastic programs with linear recourse, where for a specific realization $\xi$ of $\bm{\xi}$,
\begin{equation}\label{obj}
F(x,\xi)= c^Tx+Q(x,\xi),\\
\end{equation}
and\begin{equation}\label{recourse}
\begin{aligned}
Q(x,\xi)=&\min_y && q(\xi)^Ty\\
&\text{s.t.} && W(\xi)y\geq h(\xi)-T(\xi)x .
\end{aligned}
\end{equation}
Here $q(\xi) \in \bR^{n_2},W(\xi) \in \bR^{m_2 \times n_2},T(\xi) \in \bR^{m_2 \times n_1}$ and $h(\xi) \in \bR^{m_2}$. 
When $X$ is polyhedral, (\ref{real})-(\ref{recourse}) is a two-stage stochastic linear program.
The model \eqref{real} with $F$ defined by \eqref{obj} is said to have {\em relatively complete recourse} if there
exists a solution to \eqref{recourse} for every $x \in X$ for $\bP$-almost every $\xi\in\Xi$. Our interest in this paper is studying
this problem in the case when \eqref{real} does \emph{not} have relatively complete recourse. When \eqref{recourse} does
not have a solution, we adopt the convention that $F(x,\xi) = Q(x,\xi) = +\infty$.
Throughout this paper we assume that there exists a $\bP$-integrable random variable $Z$ such that $F(x,\bm{\xi})
\geq Z$ $\pxi$-almost surely for all $x \in X$, which assures that $f(x)$ is well-defined (but does not exclude the possibility that $f(x) = +\infty$).


Problem (\ref{real}) is usually intractable unless $\Xi$ is a small finite set
\cite{dyer2006computational,shapiro2005complexity}. A popular approach
to obtain a tractable approximation is to solve a sample average approximation (SAA) problem. The basic idea of SAA is
that, instead of solving problem (\ref{real}), the objective function is replaced by a sample average function
$N^{-1}\sum_{j \in [N]} F(x,\xi^j)$,
where we use the notation $[N] := \{1,\ldots,N\}$ for any integer $N$.
We assume throughout the paper that $\{\xi^1,\ldots,\xi^N\}$ is a sample of $N$ independent and identically distributed (iid) realizations of $\bm{\xi}$.  We call the following problem the SAA problem of (\ref{real}):
\begin{equation}\label{saadef}
\min_{x\in X}\hat{f}_N(x):=N^{-1}\sum_{j \in [N]}F(x,\xi^j)
\end{equation}

SAA has been widely studied in the contexts of convex stochastic programming \cite{king1991epi}, general stochastic
programming \cite{dupacova1988asymptotic}, stochastic discrete programming \cite{kleywegt2002sample} and two-stage
stochastic programming \cite{shapiro1998simulation}. 
The majority of these assume a finite objective, i.e., $f(x)<\infty$ for all $x\in X$, especially those related to
convergence rates or sample size estimates. In the context of two-stage stochastic programs, a necessary condition for $f(x)<\infty$ for all
$x\in X$ is that relatively complete recourse holds. 
Relatively complete recourse is often a natural assumption, in particular, since a solution that has no feasible
recourse may be considered ill-defined. From a modeling perspective, it may be argued that no matter what action is
taken in the first-stage, there should always be some feasible (potentially costly) recourse action. Indeed, it is
always possible to create a model with relatively complete recourse by
introducing auxiliary variables in the second stage that measure constraint violation and, for example, penalizing the use of
these variables in the cost function. On the other hand,
there may be recourse actions that are undesirable to use except in very rare circumstances or which are difficult to
model or solve (e.g., involving discrete decisions). In addition, penalizing constraint violation may require a high
cost of violation to assure high likelihood of feasibility, which in turn could lead to an ill-conditioned problem. 
In such a situation, a model which allows solutions that do not have a recourse action
in every possible outcome may still be meaningful, although it would be desirable to obtain a solution $x$ such that the probability of having a
recourse action, i.e., $\phi(x) := \pxi( F(x,\bm{\xi}) < \infty )$, is high. We refer to $\phi(x)$ as the {\it recourse likelihood} of a solution $x$. 

The idea of using SAA to obtain a solution that has high recourse likelihood is closely related to the scenario
approximation approach for chance constraints \cite{calafiore2005uncertain,calafiore2006scenario}. In particular, the scenario approximation proposed in \cite{calafiore2006scenario} is the problem
\begin{equation}
\label{eq:scenapp}
\min \{ h(x) : G(x,\xi^j) \leq 0, j \in [N], x \in X \}. 
\end{equation}
It is shown in \cite{calafiore2006scenario} that the measure of the set of constraints that are violated by  the
solution converges rapidly to zero. Similar results are obtained in \cite{luedtke2008sample}
without the convexity assumptions. (See Section \ref{subsect:ccp} for a detailed review of these results). If one considers $G(x,\xi)$ to be a
measure of violation of the recourse problem, then these results \emph{almost} directly apply to provide estimates of
recourse likelihood for an SAA solution from a two-stage stochastic programming problem. The challenge, however, is that
in \eqref{eq:scenapp} the function $h(\cdot)$ is deterministic, whereas in an SAA approximation of \eqref{real}, the
objective function is an expected value which is also approximated by the  sample.

In this paper we study SAA for problem \eqref{real} in the case that relatively complete recourse does not
hold, and investigate the quality of solutions obtained in terms of recourse likelihood of the solution obtained,
in addition to the expected cost. Inspired by the results on scenario approximation for chance constraints,  we first
investigate  bounds on the sample size required to obtain (with high confidence) a solution that has high recourse
likelihood (i.e., a solution $x$ with $\phi(x) \geq 1-\epsilon$) for two-stage stochastic linear programming problems. We also investigate the
use of SAA to obtain a solution with $\phi(x) = 1$, which we refer to as a {\it completely reliable} solution. When the feasible region $X$ is finite (e.g., as in a bounded pure integer program), we establish
bounds on the probability that {\it every} (near) optimal solution of the SAA problem is feasible and near-optimal to 
the true problem (these results extend similar results \cite{kleywegt2002sample,shapiro2009lectures}). For 
the more general case that $X$ is not finite, we consider two cases where a   modified SAA problem can, with high
confidence (with respect to the sample defining the SAA problem), yield a
completely reliable solution. In the first case we assume the support of the random vector is the product of the marginal supports, and in the second case we consider a two-stage stochastic linear
program in which only the right-hand side is random. Finally, we perform a numerical illustration of
the use of SAA to obtain solutions that have high recourse likelihood and padded SAA to obtain completely reliable solutions. Our numerical study confirms the viability of
these approaches, but also suggests that our sample size estimate for two-stage stochastic linear programs is not tight.
In particular, our theoretical results suggest the required sample size is $O(n_1n_2)$, whereas the numerical
study suggests the dependence on $n_2$ may not be necessary.

Aside from a few classic results on convergence in the limit and expected bias of the SAA objective value, which we
review in Section \ref{litrev}, the work that is most closely related to our results is the recent paper by Liu \cite{liu2020feasibility}, 
who establishes similar results on generating solutions with high recourse likelihood for stochastic programs
having a property they refer to as {\it chain-constrained domain}. 
In their most general results, they establish that if the problem has a chained-constrained domain ``of order $m$'', then
for $\epsilon \in (0,1)$, the probability (over the sample) that the recourse likelihood of the SAA solution is less
than $1-\epsilon$ decreases exponentially fast, provided the sample size is at least as large as $m/\epsilon$.
The dependence on $m$ may be a limitation, as for a two-stage stochastic linear program with randomness only in the
right-hand side of the constraints, $m$ is at least as large as the number of extreme rays of the subproblem feasibility
cone. Significantly stronger  results are obtained when it is assumed that $X$ is a convex set and $F(\cdot,\xi)$ is a
convex function for all $\xi \in \Xi$. Under the additional assumption that the
set of optimal solutions lies in the interior of the domain of $f(\cdot)$,
they show the very
strong result that the probability that $\phi(x_N) < 1$ converges to zero exponentially fast with $N$.
Without this assumption, they demonstrate that in the chain-constrained case the dependence on the parameter $m$ in
the convergence rate can be replaced by the number of active constraints in an optimal solution, which, e.g., can be
bounded by the number of first-stage variables $n_1$. Our results complement those in \cite{liu2020feasibility} by
conducting a different analysis which does not use the chain-constrained domain assumption nor an assumption that the set of optimal solutions lies in the interior of 
the domain of $f(\cdot)$. 
Also related, Chen et al. \cite{chen2019convergence} study SAA for two-stage stochastic generalized
equations, which can be viewed as an extension of two-stage stochastic convex programs, without assuming relatively
complete recourse. Their analysis focus on the convergence of SAA solutions to optimal solutions rather than the
recourse likelihood. 

This paper is organized as follows. In Section \ref{litrev}, we review additional related literature. In Section
\ref{saasec}, we  analyze the impact of sample size on the recourse likelihood for two-stage stochastic linear programs. In Section \ref{finitex}, we study the
rate of convergence of the set of optimal solutions of the SAA problem to the true problem in the case that the set $X$
is finite. In Section \ref{modsaa}, we introduce modified SAA programs and study its convergence properties. In Section \ref{numerical}, we present some numerical tests on a two-stage recourse planning problem.

\section{Review of related results}\label{litrev}

We discuss some existing results which are closely related to the use of SAA for obtaining approximate solutions to
\eqref{real}. We review several facts about the consistency and convergence of the SAA approach and details on related
results for chance-constrained stochastic programs.

We denote the optimal values of (\ref{real}) and (\ref{saadef}) by $v^*$ and $\hat{v}_N$, respectively. 

\subsection{SAA for two-stage stochastic programs}
The consistency of $\hat{v}_N$ and SAA solutions is proved under different assumptions in \cite{dupacova1988asymptotic,king1991epi,shapiro2009lectures}. 
For example, Theorem 5.4 of \cite{shapiro2009lectures} provides general condition under which $\hat{v}_N$ converges
to $v^*$ and SAA solutions converge to optimal solutions of \eqref{real}. This
result allows the objective $F$ to take values in $\bR\cup\{\pm\infty\}$, which includes some of the cases considered in this paper.

For the purpose of estimating the solution quality, one may be interested in the relationship between $\hat{v}^N$ and
$v^*$. For stochastic programs with a real-valued objective, it is known that $\bE[\hat{v}_N]$ is a lower bound of $v^*$
and is monotonically increasing in $N$ \cite{mak1999monte,norkin1998branch}. We do not pursue this here due to space
limitations, but we conjecture that the proof in \cite{mak1999monte} can be extended to cases where $F$ is extended real
valued provided the appropriate technical assumptions are in place to assure the compared values are well-defined.

\subsection{Scenario approximation of chance-constrained problems}\label{subsect:ccp}
The feasibility of SAA solutions are closely related to a type of optimization problems, called chance-constrained problems:\begin{equation}\label{cc}
\begin{aligned}
\min_{x\in X}\  &h(x) \\
\text{s.t.}\  &\pxi\big( G(x,\bm{\xi})\leq 0\big)\geq 1-\epsilon.
\end{aligned}
\end{equation}
Here $h:\bR^{n_1} \rightarrow \bR$ is a real-valued function on $X$ and $G(\cdot,\xi)$ is a real-valued function on $X$ for
any $\xi\in\Xi$. When $X\subseteq\bR^{n_1}$ is a closed convex set, $h(\cdot)$ is a convex function, and $G(\cdot,\xi)$
is a convex
function for any $\xi\in\Xi$, we call (\ref{cc}) a convex chance-constrained problem.

An idea similar to SAA called the scenario approach is applied in order to approximate the original convex chance-constrained problem (\ref{cc}):
\begin{equation}\label{scenappr}
\begin{aligned}
\min_{x\in X}\  &h(x) \\
\text{s.t.}\  &G(x,\xi^j)\leq 0, j \in [N].
\end{aligned}
\end{equation}
Sample size estimates for the scenario approach are studied in \cite{calafiore2006scenario,campi2008exact}. The main result of \cite{calafiore2006scenario} is as follows. Let $X^N$ denote the feasible region of (\ref{scenappr}).
\begin{thm}\label{ccresult}
	Fix two parameters $\epsilon\in(0,1)$ and $\beta\in(0,1)$. If the optimal solution of (\ref{scenappr}) is unique, and the sample size $N$ satisfies\begin{displaymath}
	N\geq\frac{2}{\epsilon}\log\frac{1}{\beta}+2(n_1+1)+\frac{2(n_1+1)}{\epsilon}\log\frac{2}{\epsilon},
	\end{displaymath}
	then the unique optimal solution $\hat{x}_N$ of (\ref{scenappr}) satisfies\begin{equation}\label{smplb}
	\pxi \big(G(\hat{x}_N,\bm{\xi})\leq 0\big)\geq 1-\epsilon
	\end{equation}
	with probability (over the sample measure $\pSamp$) at least $1-\beta$.
\end{thm}

Similar results are obtained in \cite{luedtke2008sample} without the convexity assumptions. As an example, the following
result is obtained when $X$ is finite.
\begin{thm}
	Suppose $X$ is finite and define $X_\epsilon:=\{x\in X: \pxi(G(x,\bm{\xi}) \leq 0) \geq 1-\epsilon \}.$	Then
	\begin{displaymath}
	\pSamp\big(X^N\subseteq X_\epsilon\big)\geq 1-|X\backslash X_\epsilon|(1-\epsilon)^N.
	\end{displaymath}
\end{thm}
Thus, to obtain confidence $1-\beta$ that $X^N\subseteq X_\epsilon$, we should take\begin{displaymath}
N\geq \frac{1}{\epsilon}\log\frac{1}{\beta}+\frac{1}{\epsilon}\log|X\backslash X_\epsilon|.
\end{displaymath}

\section{Recourse likelihood of SAA solutions for two-stage stochastic linear programs}\label{saasec}

In this section we consider a two-stage stochastic linear program \eqref{real} with $F(x,\xi)$ defined in \eqref{obj}-\eqref{recourse}, and 
polyhedral $X=\{x\in\bR^{n_1}:A x\leq b \}$, where $A \in\bR^{m_1\times n_1}$ and $b \in\bR^{m_1}$.
We further assume in this section that $m_2\geq n_2+1$. 
We furthermore partition the constraints in \eqref{recourse} into those that depend on the random vector $\xi$ and those
that are deterministic:  
\begin{alignat}{2}
Q(x,\xi)=&\min_y\ &&q(\xi)^Ty \nonumber \\
&\text{s.t.} && \randw(\xi)y\geq \randh(\xi)-\randt(\xi)x \label{recourse:rand} \\ 
& && Dy \geq d - Cx. \label{recourse:det}
\end{alignat}
We assume that for every $x \in X$ there exists $y$ that satisfies \eqref{recourse:det}. Since the constraints
\eqref{recourse:det} are deterministic, this assumption can always be satisfied by augmenting the definition of $X$ with
feasibility constraints of the form $\beta^\top (d - Cx) \leq 0$ for $\beta$ satisfying $D^T\beta = 0$.

Define the feasibility function $H:\mathbb{R}^{n_1}\times\mathbb{R}^d \rightarrow \mathbb{R}$ as follows: 
\begin{equation}\label{hdef}
\begin{aligned}
H(x,\xi)=\min_{\eta,y} \ &\eta\\
\text{s.t.}\ &\eta e +\randw(\xi)y\geq \randh(\xi)-\randt(\xi)x \\
& Dy \geq d- Cx. 
\end{aligned}
\end{equation}
The function $H(\cdot,\xi)$ is convex for each $\xi \in \Xi$. 
If each entry of $W(\cdot)$, $h(\cdot)$ and $T(\cdot)$ is measurable, then $H(x,\cdot)$, as the optimal value of a linear program, is measurable \cite{heilmann1977optimal} for any $x\in X$. Note that $F(x,\xi)<+\infty$ if and only if $H(x,\xi)\leq 0$.

Since $\phi(x)=1$ for any $x$ satisfying $f(x)=\bE F(x,\bm{\xi})<+\infty$, problem (\ref{real}) has the following equivalent form:
\begin{equation}\label{real2}
\begin{aligned}
\min \ &\bE F(x,\bm{\xi})\\
\text{s.t.} \ &H(x,\xi)\leq 0\quad \ \pxi\text{-almost every } \xi \in \Xi \\
&x\in X.
\end{aligned}
\end{equation}
The idea of adding constraints $H(x,\xi)\leq 0$ in the first-stage to enforce 
$F(x,\xi)< \infty$ for $\pxi$-almost every $\xi \in \Xi$ dates back to \cite{van1969shaped} (see also \cite{birge2011introduction}). 

Since $F(x,\xi)<+\infty$ if and only if $H(x,\xi)\leq 0$, the SAA problem can be written as
\begin{equation}\label{stdsaa1}
\begin{aligned}
\min \ &N^{-1}\sum_{j \in [N]}F(x,\xi^j)\\
\text{s.t.} \ &H(x,\xi^j)\leq 0,\quad j \in [N]\\
&x\in X.
\end{aligned}
\end{equation}

We investigate bounds on the probability that the solution $x_N^*$ obtained from the SAA problem \eqref{stdsaa1} has ``high recourse likelihood'', i.e., satisfies $\phi(x_N^*)\geq 1-\epsilon$.

Before we state and prove our main result, we write the (standard) Benders linear program formulation of 
(\ref{stdsaa1}).
First, observe that for $x$ and $\xi$ such that
$H(x,\xi)\leq 0$ strong duality implies
\begin{displaymath}
\begin{aligned}
Q(x,\xi)
=&\max_{\alpha \in \extr(P(\xi))}&&\alpha^T\big(h(\xi)-T(\xi)x\big),
\end{aligned}
\end{displaymath}
where $P(\xi) :=\{\alpha \in \bR^{m_2}_+ : W(\xi)^T\alpha=q(\xi) \}$
and $\extr(P)$ denotes the set of extreme points of a polyhedron $P$.

Similarly, 
\begin{alignat*}{2}
H(x,\xi)=&& \min_{\eta,y} \ &\eta\\
&&\text{s.t.}\ &\eta e_\cset+W(\xi)y\geq h(\xi)-T(\xi)x,\\
=&& \max_{\alpha\in\extr(\bar{P}(\xi))} &\alpha^T\big(h(\xi)-T(\xi)x\big),
\end{alignat*}
where $\bar{P}(\xi):=\{\alpha \in \bR_+^{m_2} : W(\xi)^T\alpha=0,e_\cset^T\alpha=1 \}$
and $\cset$ is the set of indices of the constraints \eqref{recourse:rand}.

Thus, 
introducing a variable $\Theta_j$ to represent each function $Q(\cdot,\xi^j)$ for each $j \in [N]$ in the SAA problem
yields the explicit LP 
formulation:
\begin{alignat}{2}
\min \ &c^T x + N^{-1}\sum_{j \in [N]} \Theta_j  \label{saa:benobj} \\
\text{s.t.} \ &  \Theta_j \geq \alpha^T\big(h(\xi^j)-T(\xi^j)\big), &\quad& \alpha \in \extr(P(\xi^j)), j \in [N],
\label{saa:benopt} \\
& \alpha^T\big(h(\xi^j)-T(\xi^j)x\big) \leq 0, && \alpha \in \extr(\bar{P}(\xi^j)),   j \in [N],
\label{saa:benfea}\\
&Ax \leq b. \label{saa:bencon}
\end{alignat}
Our analysis relies on considering basic solutions of the system \eqref{saa:benopt}-\eqref{saa:bencon} for subsets of
scenarios. Thus, for $\sset \subseteq [N]$ let $\bar{B}_\sset$ be the set of basic solutions to the system
\eqref{saa:benopt}-\eqref{saa:bencon} with the index set $j \in [N]$ replaced with $j \in \sset$, and let $B_\sset$ be
the projection of $\bar{B}_\sset$ to the space of $x$.
Define $\ssets_{n_1} = \big\{ \sset \subseteq [N]: |\sset| = n_1 \big\}$.
We require the
following bound on the number of basic solutions for $\sset \in \ssets_{n_1}$.

\begin{lem}
	Let $\sset \in \ssets_{n_1}$. It holds that
	\begin{equation}
	\label{eq:bibound}
	|B_\sset|\leq|\bar{B}_\sset| \leq 
	\frac{1}{n_1!}(2n_1+m_1)^{n_1}\Bigl(\frac{m_2^{n_2+1}}{(n_2+1)!}\Bigr)^{2n_1}.
	\end{equation}
\end{lem}
\noindent {\it Proof}
First observe that for any $j \in \sset$, $|\extr(P(\xi^j))|$ is bounded above by
$\binom{m_2}{m_2-k(\xi^j)}=\binom{m_2}{k(\xi)}$, where $k(\xi^j)(\leq n_2)$ is the number of linearly independent
equations in $W(\xi^j)^T\alpha=q(\xi^j)$. Therefore, $|\extr(P(\xi^j))|\leq\frac{m_2^{k(\xi)}}{k(\xi)!}\leq
\frac{m_2^{n_2}}{n_2!}$, by our assumption that $m_2 \geq n_2 + 1$.
Similarly $|\extr(\bar{P}(\xi^j))|\leq\frac{m_2^{n_2+1}}{(n_2+1)!}$ for each $j \in \sset$.

The system \eqref{saa:benopt}-\eqref{saa:bencon} with $j \in [N]$ replaced by $j \in \sset$ has $|\sset|+n_1 = 2n_1$ decision
variables, and thus a basic solution is the unique solution of a set of $2n_1$ linearly independent constraints with inequality signs replaced by equality signs. At least one of the
constraints in \eqref{saa:benopt} must be included in this set for each  $j \in
\sset$ (otherwise the decision variable $\Theta_j$ could not be obtained from solving the system), and the number of ways to choose one constraint from each is bounded above by 
$ \Bigl(\frac{m_2^{n_2}}{n_2!}\Bigr)^{n_1}$. The remaining $n_1$ linear constraints then are chosen from the remaining
constraints, the number of which is bounded above by
\[ n_1\Bigl(\frac{m_2^{n_2}}{n_2!}\Bigr)+n_1\frac{m_2^{n_2+1}}{(n_2+1)!}+m_1 \leq (2n_1+m_1)\frac{m_2^{n_2+1}}{(n_2+1)!}. \] 
Thus, the overall number of basic solutions is bounded above by
\[ \Bigl(\frac{m_2^{n_2}}{n_2!}\Bigr)^{n_1}\frac{1}{n_1!}
\Bigl((2n_1+m_1)\frac{m_2^{n_2+1}}{(n_2+1)!}\Bigr)^{n_1}\leq\frac{1}{n_1!}(2n_1+m_1)^{n_1}\Bigl(\frac{m_2^{n_2+1}}{(n_2+1)!}\Bigr)^{2n_1}.
\]
\qed

\begin{defn}
	We say $x_N^*$ is a basic optimal solution of the SAA \eqref{stdsaa1} for the two-stage stochastic LP if
	$(x,\Theta)=\Bigl(x_N^*,\bigl(Q(x_N^*,\xi^j)\bigr)_{j \in [N]}\Bigr)$ is a basic optimal solution of \eqref{saa:benobj}-\eqref{saa:bencon}.
\end{defn}

\begin{thm}\label{feathm}
	Assume the SAA (\ref{stdsaa1})
	has a basic optimal solution $x^*_N$. Then,
	\begin{displaymath}
	\pSamp( \phi(x^*_N)\geq 1-\epsilon) \geq 1- \binom{N}{n_1} \frac{1}{n_1!}(2n_1+m_1)^{n_1}\Bigl(\frac{m_2^{n_2+1}}{(n_2+1)!}\Bigr)^{2n_1}
	(1-\epsilon)^{N-n_1}.
	\end{displaymath}
\end{thm}
\noindent {\it Proof}
For each $\sset  \in \ssets_{n_1}$ define the set of samples for which $x_N^* \in B_\sset$:
\[ \Xi_\sset^N := \{ (\xi^1,\ldots,\xi^N)\in\Xi^N: x_N^* \in B_\sset  \}. \]
We first argue that
\begin{equation}
\label{eq:coverage}
\Xi^N =  \bigcup_{\sset \in \ssets_{n_1}} \Xi_\sset^N. 
\end{equation}
(This idea of partitioning the sample space is inspired by \cite{calafiore2006scenario}.)
For any fixed $(\xi^1,\ldots,\xi^N)\in\Xi^N$, $(x,\Theta)=\bigl(x_N^*,\bigl(Q(x_N^*,\xi^j)\bigr)_{j \in [N]}\bigr)$ is a
basic solution of the system \eqref{saa:benopt}-\eqref{saa:bencon}. Therefore,
$(x,\Theta)=\bigl(x_N^*,\bigl(Q(x_N^*,\xi^j)\bigr)_{j \in [N]}\bigr)$ is the unique solution of the linear system determined
by $N+n_1$ linearly independent constraints selected from \eqref{saa:benopt}-\eqref{saa:bencon}, which we refer to as
active constraints. We say scenario $j$ is a \emph{critical scenario} (see \cite{ruszczynski1997accelerating}) if at
least two of the constraints in \eqref{saa:benopt} and \eqref{saa:benfea} are active for scenario $j$. Note that at
least one of the constraints in \eqref{saa:benopt} must be active for each $j\in [N]$ (otherwise some $\Theta_j$ cannot
be determined). Therefore, there can be at most $n_1$ critical scenarios. If scenario $j$ is noncritical, variable
$\Theta_j$ can be projected out from the linear system by removing the equation that corresponds to the unique active constraint for scenario $j$. Define\begin{displaymath}
\sset:=[N]\setminus\{j_1,\ldots,j_{N-n_1} \}\in\ssets_{n_1},
\end{displaymath}
where $j_1,\ldots,j_{N-n_1}$ are $N-n_1$ noncritical scenarios. (If there are more than $N-n_1$ noncritical scenarios,
an arbitrary subset of size $N-n_1$ can be selected.) Then $\bigl(x_N^*,\bigl(Q(x_N^*,\xi^j)\bigr)_{j\in \sset}\bigr)$ is
the unique solution of the remaining linear system with $\{\Theta_j\}_{j\notin \sset}$ projected out. Therefore, by
definition, $x_N^*\in B_\sset$ for $\sset\in\ssets_{n_1}$, i.e., \eqref{eq:coverage} holds.

For any fixed $\sset \in \ssets_{n_1}$ and $x_b \in B_\sset$, we have
\begin{align*}
\pSamp\big(\phi(x_b) <&~1-\epsilon,x_N^*=x_b\big) \\
&= \pSamp\big(x_N^*=x_b\big|\phi(x_b)<1-\epsilon\big)\cdot\pSamp\big(\phi(x_b)<1-\epsilon\big)\\
&\leq\pSamp\big(H(x_b,\xi^j)\leq 0,j\notin \sset \big|\phi(x_b)<1-\epsilon\big)\cdot 1\\
&\leq(1-\epsilon)^{N-n_1}.
\end{align*}
The first inequality follows because $x_b$ satisfying the constraints $H(x_b,\xi^j) \leq 0$ for $j \notin \sset$ is a necessary condition for
$x_b$ to be a feasible solution of \eqref{stdsaa1}.
The last inequality follows because $x_b$ is the unique solution of the linear system formed by equations defined by
$\{\xi^j \}_{j\in \sset}$, which is statistically independent of $\{\xi^j\}_{j\notin \sset}$ since $\{\xi^j\}_{j \in [N]}$ is an iid sample. 
It follows from \eqref{eq:coverage} that
\begin{align*}
\pSamp\big(\phi(x_N^*)<1-\epsilon\big) & \leq \sum_{\sset \in \ssets_{n_1}} \sum_{x_b \in B_\sset}\pSamp\big(\phi(x_b)< 1-\epsilon,x_N^*=x_b\big)\\
& \leq \binom{N}{n_1}\frac{1}{n_1!}(2n_1+m_1)^{n_1}\Bigl(\frac{m_2^{n_2+1}}{(n_2+1)!}\Bigr)^{2n_1}
(1-\epsilon)^{N-n_1}.
\end{align*} \qed

\begin{cor}\label{Nbound}
	Let $\epsilon\in(0,1)$ and $\beta\in(0,1)$,
	and let $x^*_N$ be a basic optimal solution of the SAA of (\ref{stdsaa1}). 	
	If the sample size $N$ satisfies\begin{align*}
	N\geq\frac{2}{\epsilon}\log\frac{1}{\beta}&+\frac{4n_1n_2}{\epsilon}\Big(\log\bigl(\frac{m_2}{n_2+1}\bigr)+1\Big)\\
	&+\frac{2n_1}{\epsilon}\Big(\log\big(\frac{m_1}{n_1}+2\big)+\log \frac{2}{\epsilon}+1\Big)+2n_1,
	\end{align*}
	then $\pSamp(\phi(x_N^*)\geq 1-\epsilon) \geq 1-\beta$.
\end{cor}
\noindent {\it Proof} See Appendix.
\qed

As is typical in SAA analysis, Theorem \ref{feathm} is based on a number of conservative approximations, such that the
bound estimate in Corollary \ref{Nbound} is conservitive, and is hence not recommended to be used directly. Instead,
the importance of Corollary \ref{Nbound} is to demonstrate qualitatively how the sample size impacts recourse
likelihood, in particular bounding the dependence on problem size parameters such as $n_1$, $n_2$, etc. 

\noindent\textbf{Remark }When $\Xi$ is finite, define $\epsilon^*:=\min\{\bP(\bm{\xi}=\xi):\xi\in\Xi \}.$ Then, for any
$\epsilon<\epsilon^*$,  
$\phi(x_N^*)\geq 1-\epsilon$ if and only if $\phi(x_N^*)=1$. Therefore, using, e.g., $\epsilon = \epsilon^*/2$ in
Corollary \ref{Nbound} yields an estimate of the sample size required to obtain a {\it completely reliable} solution
with high confidence with respect to the sample distribution $\pSamp$. Unfortunately, this estimate is at least $\Omega(1/\epsilon^*)$, and this cannot be improved in general. 
For example, consider the following problem:
\begin{displaymath}\label{cebeg}
\min_{x\in X}f(x)=\bE F(x,\bm{\xi}),
\end{displaymath}
where $X=[0,2]$, $n\bm{\xi}$ follows the binomial distribution $B(n,1/2)$ (with support size $n+1$) and
\begin{equation}\label{eg:sample_size}
F(x,\xi)= \min \{ y : \xi \leq x\leq y \}=\left\{
\begin{array}{ll}
x, & \text{if }x\geq\xi,\\
+\infty, & \text{otherwise.}
\end{array}
\right.
\end{equation}
In this case, the solution of the approximate problem (\ref{stdsaa1}) is
\begin{displaymath}\label{ceend}
\hat{x}_N=\max_{j\in[N]}{\xi^j}.
\end{displaymath}
Then $\bP^N(f(\hat{x}_N)=+\infty)=\bP^N(\max_{j\in[N]}{\xi^j}<1)=(1-1/2^n)^N=(1-\epsilon^*)^N\approx e^{-\epsilon^*N}$,
in which case we need approximately $N\geq\log(1/\beta)/\epsilon^*=\log(1/\beta)2^n$ to have
$\bP^N(f(\hat{x}_N)<+\infty)\geq 1-\beta$. This indicates that required sample size could even grow exponentially with the
size of the support of $\Xi$.

\section{Exact convergence for finite $X$}\label{finitex}

In this section we assume the set $X$ is finite, such as in the case that all decision variables are integer and
bounded. In \cite{kleywegt2002sample}, it has been shown that, in the case when $f(x)$ is real-valued for all $x\in X$,
under certain assumptions, the optimal solution set of (\ref{stdsaa1}) converges exponentially to the optimal solution
set of (\ref{real2}). We extend this result to the extended-real-valued objective case. In particular, we obtain bounds
on the sample size required to obtain, with high confidence with respect to the sample distribution $\pSamp$, a solution that is completely reliable and nearly optimal.

We introduce notations that will be used. For $x\in X$, recall that we define $\phi(x):=\pxi(F(x,\bm{\xi})<+\infty)$ to be the recourse likelihood of $x$.
We define the set of completely reliable solutions, i.e., those which have a feasible recourse decision with
probability $1$, as 
\begin{equation}\label{xfea}
X^{\text{Fea}}:=\{x\in X:\phi(x)=1 \}
\end{equation}
and its complementary set $X^{\text{Infea}}:=\{x\in X:\phi(x)<1 \}$.
By definition, the effective domain of $f(\cdot)$ is a subset of $X^{\text{Fea}}$. For any $\epsilon\geq0$, we denote the sets of $\epsilon$-optimal solutions of (\ref{real}) and (\ref{saadef}) by $S^\epsilon$ and $\hat{S}^\epsilon_N$, respectively, i.e.
\begin{displaymath}
\begin{aligned}
S^\epsilon:=\{x \in X :f(x)\leq v^*+\epsilon \}, \quad
\hat{S}^\epsilon_N:=\{x \in X :\hat{f}_N(x)\leq \hat{v}_N+\epsilon\}.
\end{aligned}
\end{displaymath}
Note that when $f(x)=+\infty$ for all $x\in X$, i.e., $v^*=+\infty$, this definition implies $S^\epsilon=X$ for any
$\epsilon \geq 0$. Under this convention, when $v^*=+\infty$ the inclusion
$\hat{S}^\delta \subseteq S^\epsilon$ holds trivially for any $\delta \in [0,\epsilon]$, and hence the results we
present in this section are trivially correct in that case.

Our results on convergence of the set of optimal solutions $\hat{S}^\delta_N$ and optimal value $\hat{v}_N$ require large
deviations (LD) theory. We first review some of the LD theory following the presentation in \cite{kleywegt2002sample} and \cite{shapiro2009lectures}.
Consider a random variable $Y$ with mean $\mu=\bE[Y]$. Let $Y_1,\ldots,Y_N$ be an iid sequence of $N$ realizations of
the random variable $Y$ and the average $Z_N:=N^{-1}\sum_{j \in [N]} Y_j$. Let $M(t):=\bE[e^{tY}]$ denote the moment-generating function of $Y$. The following inequality holds:
\begin{equation}\label{Izineq}
\pSamp(Z_N\geq a)\leq e^{-NI(a)},
\end{equation}
where $I(z):=\sup_t\{tz-\log[M(t)]\}$ is the conjugate of the logarithmic moment-generating function.
Furthermore, from \cite{dembo2010large}, if we assume the moment-generating function $M(t)$ is finite valued in a neighborhood of 0, it follows that $I(a)>0$ for all $a\neq\mu$.

Next we implement this LD theory to the analysis of our problem.
First, let $u:X^{\text{Fea}}\backslash S^\epsilon\rightarrow X^{\text{Fea}}$ be a mapping such that for some $\epsilon^*>\epsilon$,
\begin{displaymath}
\begin{aligned}
f(u(x))\leq f(x)-\epsilon^* \qquad\forall x\in X^{\text{Fea}}\backslash S^\epsilon.
\end{aligned}
\end{displaymath}
Since $X^{\text{Fea}}$ is finite, such $u(\cdot)$ exists. One example is $u:X^{\text{Fea}}\backslash S^\epsilon\rightarrow S^0$ with $\epsilon^*=\min_{x\in X^{\text{Fea}}\backslash S^\epsilon} f(x)-v^*>\epsilon$.
The following result follows from \cite{kleywegt2002sample}.
\begin{lem}
	Let $\epsilon$ and $\delta$ be nonnegative numbers such that $\delta\in[0,\epsilon]$. Then for any $x\in
	X^{\text{Fea}}\backslash S^\epsilon$, $\pSamp(x\in\hat{S}_N^\delta )\leq e^{-NI_x(-\delta)}$ where
	\begin{displaymath}
	I_x(z):=\sup_t\big\{tz-\log\bE\big[e^{t(F(x,\bm{\xi})-F(u(x),\bm{\xi})}\big] \big\}.
	\end{displaymath}
\end{lem}

We make the following assumption:

\begin{assumption}
	\label{assum:mgf}
	For every $x\in X^{\text{Fea}}\backslash S^\epsilon$, the moment-generating function of the random variable $F\big(u(x),\bm{\xi}\big)-F(x,\bm{\xi})$ is finite valued in a neighborhood of 0. 
\end{assumption}

Assumption \ref{assum:mgf} is equivalent to the distribution of $F\big(u(x),\bm{\xi}\big)-F(x,\bm{\xi})$ being
sub-exponential (see \cite{vershynin2018high}). For example, Assumption \ref{assum:mgf} holds when $\Xi$ is bounded.
Later we make the stronger Assumption \ref{assum:strmgf} for the special case when the distribution of $F\big(u(x),\bm{\xi}\big)-F(x,\bm{\xi})$ is sub-Gaussian.

Let $\eta:=\min\{\pxi(F(x,\bm{\xi})=\infty):x\in X^{\text{Infea}}\}$, then we have the following lemma which bounds the
likelihood that a solution $x$ with $\phi(x)<1$ is feasible to the SAA problem.
\begin{lem}\label{lem1}
	For every $x\in X^{\text{Infea}}$, $\pSamp\big(\hat{f}_N(x)<+\infty\big)\leq(1-\eta)^N$.
\end{lem}
\noindent {\it Proof}  By definition of $\eta$, $\pxi\big(F(x,\bm{\xi})=\infty\big)\geq\eta$. So for $j \in [N]$, we have $\pxi\big(F(x,\xi^j)<\infty\big)\leq 1-\eta$. Since $\xi^1,\ldots,\xi^N$ are independent,
\begin{displaymath}
\begin{aligned}
\pSamp\big(\hat{f}_N(x)<+\infty\big)
=&\pSamp\Big(\bigcap_{j \in [N]}\{F(x,\xi^j)<\infty\}\Big)\\
=&\prod_{j \in [N]}\pxi\big(F(x,\xi^j)<\infty\big)
\leq(1-\eta)^N.
\end{aligned}
\end{displaymath}
\qed

We next show that the probability that a $\delta$-optimal solution to the SAA problem is not an $\epsilon$-optimal
solution to problem decreases to zero exponentially fast.

\begin{thm}
	Let $\epsilon>0$ and $\delta\in[0,\epsilon]$. Then
	\begin{displaymath}
	\begin{aligned}
	\pSamp\big(\hat{S}^\delta_N\nsubseteq S^\epsilon\big) 
	&\leq |X^{\text{Infea}}|e^{-N\eta}+|X^{\text{Fea}}\backslash S^{\epsilon}|e^{-N\gamma(\delta,\epsilon)},
	\end{aligned}
	\end{displaymath}
	where
	\begin{displaymath}
	\gamma(\delta,\epsilon)=\min_{x\in X^{\text{Fea}}\backslash S^{\epsilon}}I_x(-\delta).
	\end{displaymath}
	Moreover, if Assumption \ref{assum:mgf} holds, then $\gamma(\delta,\epsilon)>0$.
\end{thm}
\noindent {\it Proof}  For $x\in X^{\text{Infea}}$, by Lemma \ref{lem1}
\begin{displaymath}
\pSamp\big(x\in \hat{S}^\delta_N\big)\leq\pSamp\big(x\text{ is a feasible solution of (\ref{stdsaa1})}\big)\leq(1-\eta)^N.
\end{displaymath}
For $x\in X^{\text{Fea}}\backslash S^\epsilon$, $\pSamp\big(x\in \hat{S}^\delta_N\big)\leq e^{-NI_x(-\delta)}$.
Therefore,
\begin{displaymath}
\begin{aligned}
\pSamp\big(\hat{S}^\delta_N\nsubseteq S^\epsilon\big)
\leq&\sum_{x\in X\backslash S^\epsilon}\pSamp\big(x\in \hat{S}^\delta_N\big)\\
=&\sum_{x\in X^{\text{Infea}}}\pSamp\big(x\in \hat{S}^\delta_N\big)+\sum_{x\in X^{\text{Fea}}\backslash S^\epsilon}\pSamp\big(x\in \hat{S}^\delta_N\big)\\
\leq&|X^{\text{Infea}}|(1-\eta)^N+|X^{\text{Fea}}\backslash S^{\epsilon}|e^{-N\gamma(\delta,\epsilon)}\\
\leq&|X^{\text{Infea}}|e^{-N\eta}+|X^{\text{Fea}}\backslash S^{\epsilon}|e^{-N\gamma(\delta,\epsilon)}.
\end{aligned}
\end{displaymath}
Under Assumption \ref{assum:mgf}, since $\delta<\bE[F(x,\bm{\xi})-F(u(x),\bm{\xi})]$ and $X^{Fea}\setminus S^{\epsilon}$ is finite, we have
\begin{displaymath}
\gamma(\delta,\epsilon)=\min_{x\in X^{\text{Fea}}\backslash S^{\epsilon}}I_x(-\delta)>0.
\end{displaymath}
\qed

Let $\tilde{\gamma}(\delta,\epsilon)=\min\{\eta,\gamma(\delta,\epsilon)\}$, then
\begin{displaymath}
\pSamp\big(\hat{S}^\delta_N\nsubseteq S^\epsilon\big)\leq|X\backslash S^\epsilon|e^{-N\tilde{\gamma}(\delta,\epsilon)}.
\end{displaymath}
For any $\beta\in(0,1)$, with sample size
\begin{displaymath}
N\geq\frac{1}{\tilde{\gamma}(\delta,\epsilon)}\log\frac{|X\backslash S^\epsilon|}{\beta},
\end{displaymath}
we have $\pSamp\big(\hat{S}^\delta_N\subseteq S^\epsilon\big)\geq 1-\beta$.

A stronger version of Assumption \ref{assum:mgf} can lead to an explicit lower bound for $\gamma(\delta,\epsilon)$:
\begin{assumption}
	\label{assum:strmgf}
	There exists a constant $\sigma>0$ such that for every $x\in X^{\text{Fea}}\backslash S^\epsilon$, the moment-generating function $M_x(\cdot)$ of the random variable $F\big(u(x),\bm{\xi}\big)-F(x,\bm{\xi})-\bE\big[F\big(u(x),\bm{\xi}\big)-F(x,\bm{\xi})\big]$ satisfies $M_x(t)\leq e^{\sigma^2t^2/2}$ for all $t\in\bR$.
\end{assumption}

If Assumption \ref{assum:mgf} is replaced by Assumption \ref{assum:strmgf}, we have $\gamma(\delta,\epsilon)>\frac{(\epsilon-\delta)^2}{2\sigma^2}$ (see \cite[page 183]{shapiro2009lectures}).

\begin{thm}
	Suppose Assumption \ref{assum:strmgf} holds. Let $\epsilon$ and $\delta$ be nonnegative numbers such that $\delta<\epsilon$. Then
	\begin{displaymath}
	\pSamp(\hat{S}^\delta_N\nsubseteq S^\epsilon)\leq|X^{\text{Infea}}|e^{-N\eta}+|X^{\text{Fea}}\backslash S^{\epsilon}|e^{-\frac{N(\epsilon-\delta)^2}{2\sigma^2}}.
	\end{displaymath}
\end{thm}

Note that results in this section do not depend on whether the two-stage stochastic program has linear recourse, and therefore can be applied to more general stochastic programs.

\section{Using SAA to obtain completely reliable solutions}\label{modsaa}

In Section \ref{finitex}, we showed that the probability (with respect to the sample distribution $\pSamp$) the SAA yields a completely reliable solution approaches
one exponentially fast, in the case that the set $X$ is finite. Unfortunately, such a result is not possible in general.
Even in the case when $\Xi$ is finite, we have seen at the end of Section \ref{saasec} that we may need a sample size proportional to $1/\epsilon^*$ with $\epsilon^*=\min\{\bP(\bm{\xi}=\xi):\xi\in\Xi \}$, which is basically saying we will need to sample almost all points of $\Xi$. When $|\Xi|=+\infty$, we may have no chance of obtaining a completely reliable SAA solution for any finite $N$. For example, consider the example \eqref{eg:sample_size} with $\bm{\xi}=\frac{\bm{\xi'}}{1+\bm{\xi'}}$ and $\bm{\xi'}$ following a Poisson distribution. We have $\pSamp(f(\hat{x}_N)=+\infty)=\pSamp(\max_{j\in[N]}{\xi^j}<1)=1$ for any finite $N$. We thus explore in this section how modified
SAA problems can be used to obtain a completely reliable solution with high confidence in some cases.   

\sloppypar{
	Recall the formulation \eqref{real2}, where $H(x,\xi)$ is defined in \eqref{hdef}.
	To simplify analysis, we
	consider an 
	alternative formulation which instead enforces ${H(x,\xi)\leq 0}$ for {\it all} $\xi \in \Xi$:
	\begin{equation}\label{real3}
	\begin{aligned}
	&\min_{x \in X^*} && \bE F(x,\bm{\xi}),
	\end{aligned}
	\end{equation}
	where $X^*:=\{x\in X:H(x,\xi)\leq 0 \text{ for all } \xi\in\Xi \}$.
	These two problems are equivalent in many cases, for example, when $H(x,\cdot)$ is lower semi-continuous.
	
	\begin{prop}
		\label{prop:hsc}
		Suppose $H(x,\cdot)$ is lower semi-continuous for all $x\in X$. Then	problem (\ref{real2}) is equivalent to problem (\ref{real3}).
	\end{prop}
	\noindent {\it Proof}  We only need to prove that the constraints are equivalent in (\ref{real2}) and (\ref{real3}). A
	feasible solution of (\ref{real3}) is trivially a feasible solution of (\ref{real2}). For $x\in X$ satisfying
	$H(x,\xi)\leq 0$ for $\pxi$-almost every $\xi \in \Xi$, let $\Xi_x=\{\xi \in \Xi : H(x,\xi)\leq 0 \}$.
	Then $\pxi(\bm{\xi}\in \Xi_x)=1$ and $\Xi_x$ is closed because of lower semi-continuity of $H(x,\xi)$. Because
	the support $\Xi$
	is the smallest closed set satisfying $\pxi(\bm{\xi}\in \Xi)=1$, we have $\Xi\subseteq \Xi_x.$ So $x$ satisfies $H(x,\xi)\leq 0$ for all $\xi\in\Xi$.
	\qed
	\noindent The condition that $H(x,\cdot)$ is lower semi-continuous is satisfied when the problem \eqref{hdef} has fixed recourse (i.e., $W(\xi)\equiv W$ is a
	fixed matrix) and both $h(\xi)$ and $T(\xi)$ are continuous in $\xi$. Indeed, in this case the dual of \eqref{hdef} is the pointwise
	maximum of finitely many continuous functions of $\xi$, which implies $H(x,\cdot)$ is continuous for all $x\in X$. See
	\cite{wets1985continuity} for more general conditions under which the optimal objective value of a linear program is lower semi-continuous in its parameters. 
	
	When the condition of Proposition \ref{prop:hsc} does not hold, we can still say that \eqref{real3} is a conservative
	approximation of of \eqref{real2}, i.e., a feasible solution of \eqref{real3} is guaranteed to be feasible to
	\eqref{real2}. Note that $x \in X^*$ implies that a recourse action exists for every possible scenario, but
	does not guarantee that $f(x)=\mathbb{E}F(x,\xi)$ is finite in general. Additional assumptions are required to ensure this.
	For example, for two-stage stochastic linear programs with fixed recourse, if $\xi$ satisfies a weak covariance condition \cite{wets1974stochastic}, then $f(x)<+\infty$ for all $x\in X^*$.
	
	Replacing the objective in \eqref{real3} with a sample average approximation, but keeping the
	feasible region yields the approximation:
	\begin{equation}\label{realsaa}
	\min_{x\in X^*}N^{-1}\sum_{j \in [N]} F(x,\xi^j) .
	\end{equation}
	Since the feasible region in \eqref{realsaa} is not approximated by sampling, standard SAA convergence results for stochastic
	programs with a real-valued objective function can be directly applied to this case.
	Directly solving (\ref{realsaa}) with $X^*$
	defined by infinite number of constraints is challenging. The tractability of such sets is discussed in robust
	optimization, e.g., \cite{margellos2014road}. In some cases, including two-stage stochastic programs with fixed
	recourse, $H(x,\cdot)$ is
	convex for any $x\in X$. In addition, if $\Xi$ is known to be a polytope, then $X^*$ can be represented by $H(x,\xi)\leq
	0$ for $\xi$ in the (finite) set of extreme points of $\Xi$.
	However, even in this case the number of constraints is potentially exponential in $d$, and separating these constraints
	requires maximizing a convex function, and is thus computationally challenging in general. Other techniques in robust
	optimization \cite{bertsimas2006tractable} or semi-infinite programs \cite{mitsos2011global} can be applied to
	approximate the feasible set $X^*$.  
	
	The focus in the remainder of this section is on the case where solving \eqref{realsaa} is computationally intractable, or if $\Xi$ is not known explicitly and instead
	we only have access to samples of $\xi$. In this case, we propose modified SAA problems
	that can yield completely reliable solutions. We study two
	cases when we can find such solutions. In section \ref{hcsec}, we consider the case when the
	support $\Xi$ of $\xi$ is equal to the product of the supports of its components. In section \ref{rhssec}, we consider stochastic linear programs in which only the right-hand side
	is random.
	
	Both cases require the following strict feasibility assumption for the problem \eqref{real3}, which we make for the
	remainder of this section.
	\begin{assumption}
		\label{assum:hfeas}
		There exists $\bar{\gamma} > 0$ such that the set $\{ x \in X : H(x,\xi) + \bar{\gamma} \leq 0, \ \forall \xi \in \Xi \}
		\neq \emptyset$.
	\end{assumption}
	Recalling the definition of $H(x,\xi)$ in \eqref{hdef}, if the constraints \eqref{recourse:rand}
	include equality constraints, satisfying this assumption would require substituting out enough decision
	variables to eliminate the equality constraints. Verifying Assumption \ref{assum:hfeas} can be difficult in
	general. An example of a sufficient condition that can be checked is if a set $\bar{\Xi}$ is known to contain the
	support $\Xi$, one can seek a feaible solution to the adjustable robust optimization problem $\inf_{x\in X}\sup_{\xi\in\bar{\Xi}}H(x,\xi)$. 
	If a solution $x \in X$ can be found with negative objective value then Assumption \ref{assum:hfeas} holds. If the problem has fixed
	recourse and $\bar{\Xi}$ is polyhedral, one approach for attempting to identify such a solution is to use affine decision rules \cite{ben2004adjustable}. 
	
	\subsection{Two-stage stochastic programs with product of marginal supports}\label{hcsec}
	\label{sec:prodmarg}
	
	We make the following assumption in this subsection.
	
	\begin{assumption}
		\label{assum:lip}
		$H(x,\cdot)$ is a Lipschitz continuous function for all $x\in X$ under infinity norm, with a uniform Lipschitz constant $L$ independent of $x$.
	\end{assumption}
	Assumption \ref{assum:lip} holds for two-stage stochastic programs under mild conditions.
	
	\begin{prop}\label{lpprop}
		Consider the two-stage stochastic program  (\ref{real})-(\ref{obj}) 
		and $H(x,\xi)$ defined by \eqref{hdef}.
		Suppose $(W(\xi),T(\xi),h(\xi))$ is Lipshitz continuous in $\xi$. Then Assumption \ref{assum:lip} holds for some $L>0$ if one of the following conditions holds:\begin{enumerate}
			\item The problem has only right-hand side randomness, i.e., $\randw(\xi)\equiv \randw$ and $\randt(\xi)\equiv
			\randt$ do not depend on $\xi$;
			\item Set $X$ is bounded and the problem has fixed recourse, i.e., $\randw(\xi)\equiv \randw$ does not depend on $\xi$;
			\item The set $\{(x,y):x\in X, Dy \geq d - Cx \}$ is bounded.
		\end{enumerate}
	\end{prop}
	\noindent {\it Proof} See Appendix.
	\qed
	
	We make the following assumption on the support of the random vector.
	
	\begin{assumption}
		\label{assum:prod}
		The support $\Xi$ of $\bm{\xi}$ is bounded and equal to the Cartesian product of supports $\Xi_i$ of $\bm{\xi}_i$,
		i.e., $\Xi=\prod_{i \in [d]}\Xi_i$.
	\end{assumption}
	
	Note that in Assumption \ref{assum:prod} we are not assuming we know $\Xi_i$ for any $i$.
	
	Now let $\xi^1,\ldots,\xi^N$ be an iid sample of $\bm{\xi}$. Given this sample, we define 
	the ``mixed scenario" $\xi^\sset=(\xi_1^{i_1},\ldots,\xi_d^{i_d})$ for each $\sset=(i_1,\ldots,i_d)\in [N]^d$. Under
	Assumption \ref{assum:prod}, if $\xi^1,\ldots,\xi^N\in\Xi$, then $\xi^\sset\in\Xi$. Therefore, for a potential solution $x$,
	it is valid to enforce $H(x,\xi^\sset)\leq 0$ for any mixed scenario $\xi^\sset$. In addition, we consider to add a $\gamma$-padding in the feasibility constraints. Then our padded SAA problem is:
	\begin{equation}\label{modapprox}
	\begin{aligned}
	&\min_x && N^{-1}\sum_{j \in [N]} F(x,\xi^j)\\
	&\text{s.t.} && H(x,\xi^\sset)+ \gamma\leq 0, &&\sset\in[N]^d,\\
	& && x\in X.
	\end{aligned}
	\end{equation}
	Let $\hat{X}_{N,\gamma}$ denote the feasible set of \eqref{modapprox}.
	When Assumption \ref{assum:prod} is not satisfied, (\ref{modapprox}) can be seen as a conservative approximation of
	problem (\ref{real3}). 
	
	\subsubsection{Solving \eqref{modapprox}}\label{subsubsec:solving_padded}
	
	In the case of two-stage linear programs \eqref{real}-\eqref{recourse}, the padded SAA problem \eqref{modapprox} can be
	reformulated as an explicit LP
	by introducing auxiliary variables to represent the recourse action in each scenario $\xi^\sset$, for $\sset \in [N]^d$: 
	\begin{equation}\label{LP:padded}
	\begin{aligned}
	&\min_{x,y,z} && c^Tx+N^{-1}\sum_{j \in [N]} q(\xi^j)^Ty(\xi^j)\\
	&\text{s.t.} && W(\xi^j)y(\xi^j)+T(\xi^j)x\geq h(\xi^j), &&j \in [N],\\
	&&&W(\xi^\sset)z(\xi^\sset)+T(\xi^\sset)x\geq h(\xi^\sset) + \gamma e_{\cset}, &&\sset \in [N]^d,\\
	& && x\in X.
	\end{aligned}
	\end{equation}
	The formulation \eqref{LP:padded} can be difficult to solve due to its exponential size. 
	However, if $H(x,\cdot)$ satisfies a monotonicity assumption, padding for one particular $\xi^\sset$ is sufficient for solving $(\ref{modapprox})$.
	
	\begin{exmp}\label{eg:monotone}
		Suppose that $H(x,\xi)$ is monotone in $\xi$, i.e., $H(x,\underline{\xi})\leq
		H(x,\bar{\xi})$ for all $x\in X$ and $\underline{\xi},\bar{\xi}\in\Xi$ satisfying $\underline{\xi}\leq\bar{\xi}$. 
		Let $\xi^{\max}$ be the vector defined by $\xi^{\max}_i = \max\{ \xi_i^j : j\in [N] \}$. 
		Then by the monotonicity assumption, $H(x,\xi^{\max})+\gamma\leq 0$ dominates $H(x,\xi^\sset)+\gamma\leq 0$ for all
		$\sset \in [N]^d$. In this case, 
		the feasible region of (\ref{modapprox}) is 
		simplified to
		\[ \hat{X}_{N,\gamma} = \{ x \in X : H(x,\xi^{\max}) + \gamma \leq 0 \} . \]  
		One example where $H(x,\xi)$ is monotone in $\xi$ is when $X \subseteq \mathbb{R}^{n_1}_+$, the constraints $y \geq
		0$ are implied by the constraints in \eqref{hdef},
		and all entries of $(-W(\xi),-T(\xi),h(\xi))$ are monotone in $\xi$. In this case, the exponential set of constraints
		in \eqref{LP:padded} can be replaced by the single set of constraints:
		\[ W(\xi^{\max})z(\xi^{\max})+T(\xi^{\max})x\geq h(\xi^{\max})+\gamma e_\cset \]
		with the single set of auxiliary variables $z(\xi^{\max})$. Note that this monotonicity assumption always holds if
		$h(\xi)=\xi$ and $W(\xi) \equiv W$ and $T(\xi) \equiv T$ are fixed.
	\end{exmp}
	
	We next discuss how \eqref{modapprox} can be solved when the monotonicity assumption in Example \ref{eg:monotone} is not
	satisfied.
	A natural approach is to replace the index set $[N]^d$ by a small subset $\mathcal{\sset}$ (could be initialized as empty)
	and iteratively add elements to the set $\mathcal{\sset}$ as needed. Specifically, after solving the problem
	\eqref{LP:padded} with index set $[N]^d$ replaced by $\mathcal{\sset}$, suppose $\hat{x}$ is the optimal
	solution. Then, we solve the following auxiliary problem 
	\begin{equation}
	\label{eq:sep}
	\hat{v} = \max \{ H(\hat{x},\xi^\sset) : \sset \in [N]^d \}
	\end{equation}
	and let $\hat{\sset}$ be an index set that attains the maximum.
	If $\hat{v} + \gamma \leq 0$, then the solution $\hat{x}$ is feasible to \eqref{LP:padded} (and hence \eqref{modapprox}),
	and hence is optimal. Otherwise, we add $\hat{\sset}$ to $\mathcal{\sset}$ and repeat the process until a feasible solution is
	obtained. The key question in this algorithm is how to solve the problem \eqref{eq:sep}. The following proposition
	provides a mixed-integer linear programming (MILP) formulation for solving this in the
	case that all the data are linear functions of $\xi$.
	
	\begin{prop}\label{prop:separation}
		Assume each entry of $\randw(\xi)$, $\randt(\xi)$ and $\randh(\xi)$ is linear in $\xi$, i.e.,\begin{enumerate}
			\item For each $k \in [n_1]$ there exists $\randt^k\in\bR^{m_2\times d}$ such that $[\randt(\xi)]_{\cdot
				k}=\randt^k\xi$;
			\item For each $k \in [n_2] $ there exists $\randw^k\in\bR^{m_2\times d}$ such that $[\randw(\xi)]_{\cdot
				k}=\randw^k\xi$;
			\item There exists $\overline{H}\in\bR^{m_2\times d}$ such that $\randh(\xi)=\overline{H}\xi$.
		\end{enumerate}
		Then the problem \eqref{eq:sep} 
		can be formulated as the following MILP:
		\begin{alignat}{3}
		&&\max_{\alpha,\beta,\delta,z}\ &\sum_{p \in I}\sum_{q \in [d]}\sum_{j \in [N]} \xi_q^j \Bigl( \overline{H}_{pq}
		-  \sum_{k \in [n_1]}\hat{x}_k \randt^k_{pq} \Bigr) z_{pqj} &&
		+ \beta^\top (d - D\hat{x}) \nonumber\\
		&&\text{s.t. }&
		\sum_{p \in I} z_{pqj} = \delta_{qj}, && q \in [d],~j \in [N] \nonumber\\
		&&&\sum_{j \in [N]} z_{pqj}  = \alpha_p, && p \in I, \ q \in [d], \nonumber\\
		&&&\sum_{p \in I} \sum_{q \in [d]} \sum_{j \in [N]} \randw^k_{pq} \xi_q^jz_{pqj} + \beta^T D^k=0,&&k \in [n_2],\label{con:old_dual}\\
		&&&e^T\alpha=1, \nonumber\\
		&&&\alpha\geq 0,~\beta\geq0,~z\geq 0,~\delta\in\{0,1\}^{d\times N}\nonumber
		\end{alignat}
		where $D^k$ is column $k$ of $D$ for $k \in [n_2]$.
		In particular, if $(\alpha^*,\beta^*,\delta^*,z^*)$ is an optimal solution of this MILP, then $\sset^*=(j^*_1,\ldots,j^*_d)$
		is an optimal solution of \eqref{eq:sep}, where for each $q \in [d]$, $j^*_q$ is the unique index in $[N]$ such that
		$\delta^*_{qj^*_q}=1$. 
	\end{prop}
	\noindent {\it Proof} See Appendix.
	\qed
	
	Although the assumptions of Proposition \ref{prop:separation} appear technical, they are
	satisfied in the common case that the components of the matrices $T$,$W$, and $h$ are mutually independent random variables
	and we simply define $\xi=(T,W,h)$. In this case the matrices $\randt^k$, $\randw^k$, and $\randh$  
	in the assumptions of Proposition \ref{prop:separation} simply represent the relevant projections of this ``flattened''
	random vector $\xi$. 
	
	When the problem has fixed recourse, i.e., $\randw(\xi)$ is deterministic, we propose a strengthened and more compact MILP formulation for problem \eqref{eq:sep}.
	\begin{prop}\label{prop:separation+}
		Assume $\randw(\xi)\equiv\randw$ is deterministic, and each entry of $\randt(\xi)$ and $\randh(\xi)$ is linear in $\xi$, i.e.,\begin{enumerate}
			\item For each $k \in [n_1]$ there exists $\randt^k\in\bR^{m_2\times d}$ such that $[\randt(\xi)]_{\cdot
				k}=\randt^k\xi$;
			\item There exists $\overline{H}\in\bR^{m_2\times d}$ such that $\randh(\xi)=\overline{H}\xi$.
		\end{enumerate}
		Define $\xi_q^{\min}:=\min_{j\in[N]}\xi_q^j$ and $\xi_q^{\max}:=\max_{j\in[N]}\xi_q^j$ for $q\in[d]$.	Then the problem \eqref{eq:sep} 
		can be formulated as the following MILP:
		\begin{align}
		&&\max_{\alpha,\beta,\delta,z,w}\ &\sum_{p\in I}\sum_{q\in[d]}\Big[\Big(\overline{H}_{pq}-\sum_{k\in[n_1]}\hat{x}_k\randt^k_{pq}\Big)(\xi_q^{\min}z_{pq1}+&&\hspace{-5.5mm}\xi_q^{\max}z_{pq2}) \Big]+ \beta^T(d - C\hat{x})\nonumber\\
		&&\text{s.t. }&\sum_{p\in I}z_{pqj}=\delta_{qj},&&\hspace{-12	mm}q\in[d],~j\in\{1,2\},\nonumber\\
		&&&z_{pq1}+z_{pq2}=\alpha_p,&&\hspace{-12mm}p\in I,~q\in[d],\nonumber\\
		&&&w_{pq1}+w_{pq2}=\beta_p,&&\hspace{-12mm}p\in [m_2]\setminus I,~q\in[d],\nonumber\\
		&&&\sum_{p\in I}\randw_{pk}z_{pqj}+\sum_{p\in [m_2]\setminus I}D_{pk}w_{pqj}=0,&&\hspace{-12mm}k\in[n_1],~q\in[d],~j\in\{1,2\},\label{con:new_dual}\\
		&&&e^T\alpha=1,\nonumber\\
		&&&\alpha\geq 0,~\beta \geq 0,~z\geq 0,~w\geq 0, ~\delta\in\{0,1\}^{d\times 2}.\nonumber
		\end{align}
		In particular, if $(\alpha^*,\beta^*,\delta^*,z^*,w^*)$ is an optimal solution of this MILP, then $\sset^*=(j^*_1,\ldots,j^*_d)$
		is an optimal solution of \eqref{eq:sep}, where for each $q \in [d]$, $j_q^*$ satisfies $\xi^{j_q^*}_q=\xi^{\min}_q$ if
		$\delta^*_{q1}=1$ or $\xi^{j_q^*}_q=\xi^{\max}_q$ if
		$\delta^*_{q2}=1$.
	\end{prop}
	\noindent {\it Proof} See Appendix.
	\qed
	
	In the case of fixed recourse, the formulation in Proposition \ref{prop:separation+} has $d\times 2$ binary variables,
	reduced from $d\times N$ in the formulation of Proposition \ref{prop:separation}.
	This reduction is made possible by exploiting convexity of $H(\hat{x},\cdot)$ which implies the maximum in
	\eqref{eq:sep} can be replaced by the maximum over the extreme points of the set $\prod_{q \in [d]}[\xi^{\min}_q, \xi^{\max}_q]$. Fixed recourse also
	enables application of the reformulation-linearization technique \cite{sherali2013reformulation}, which introduces new
	variables $\beta$ and strengthens constraints \eqref{con:old_dual} to \eqref{con:new_dual}. In Section
	\ref{subsec:testpaddedSAA_nonmon} we present results from numerical experiments that indicate the formulation
	from Proposition \ref{prop:separation+} is much more tractable than that from Proposition \ref{prop:separation},
	demonstrating the importance of exploiting the fixed recourse assumption when possible.
	
	When the monotonicity assumption in Example \ref{eg:monotone} is not satisfied but the feasibility function
	$H(\hat{x},\xi)$
	is ``partially monotone" in $\xi$, it is possible to simplify the MILP formulations for solving \eqref{eq:sep} by introducing 
	binary decision variables to determine the mixed scenario only for the components of $\xi$ for which $H(x,\xi)$ is not
	monotone. An example is presented in Section \ref{subsec:testpaddedSAA_nonmon}.
	
	Finally, we mention that since the problem size of \eqref{LP:padded} 
	can still be large even when the index set $[N]^d$ is replaced with a smaller set $\mathcal{\sset}$, one may apply standard Benders decomposition
	\cite{bnnobrs1962partitioning,van1969shaped} to solve it. In particular, 
	we can eliminate the auxiliary variables $z(\xi^J)$ and the padded feasibility constraints for ``mixed scenarios"
	$\xi^J$ for $J\in\mathcal{\sset}$ and instead enforce these constraints by adding Benders cuts as needed
	within a cutting plane algorithm.

	\subsubsection{Obtaining completely reliable solutions from \eqref{modapprox}}
	\label{sec:comprel}
	
	\begin{thm}\label{prob}
		Suppose that Assumptions \ref{assum:lip} and \ref{assum:prod} hold, and the diameter of $\Xi_i$ is $D_i$. Let $D:=\max_i D_i$ and $\eta$ be a constant defined in (\ref{etadef}). Then $\hat{X}_{N,\gamma}$  satisfies
		\begin{displaymath}\label{pr2}
		\pSamp\big(\hat{X}_{N,\gamma}\subseteq X^{\text{Fea}}\big)\geq 1-(dDL/\gamma)(1-\eta)^N.
		\end{displaymath}
	\end{thm}
	\noindent {\it Proof}  For each $i$, there exists a  $(\gamma/2L)$-net of $\Xi_i$, i.e., there exists $\bar{\xi}^k_i,\
	k\in [M_i]$, such that
	\begin{displaymath}
	\begin{aligned}
	&\bigcup_{k \in [M_i]} \big[\bar{\xi}^k_i-\gamma/2L,\bar{\xi}^k_i+\gamma/2L\big]\supseteq \Xi_i,\\
	&\pxi\big(\xi_i\in\big[\bar{\xi}^k_i-\gamma/2L,\bar{\xi}^k_i+\gamma/2L\big]\big)>0.
	\end{aligned}
	\end{displaymath}
	We can choose $\bar{\xi}^k_i,\ k\in [M_i]$ in a way that
	$M_i\leq DL/\gamma.$
	Define
	\begin{equation}\label{etadef}
	\eta(\gamma):=\min_{i,k}\pxi\big(\xi_i\in\big[\bar{\xi}^k_i-\gamma/2L,\bar{\xi}^k_i+\gamma/2L\big]\big).
	\end{equation}
	Let $B_r(\xi)$ denote the infinity norm ball of radius $r$ and center $\xi$ in $\bR^d$  and let $B^K$ denote the
	infinity norm ball $B_{\gamma/2L}\big(\big(\bar{\xi}_1^{k_1},\ldots,\bar{\xi}_d^{k_d}\big)\big)$ where
	$K=(k_1,\ldots,k_d)\in\prod_{i \in [d]}[M_i]$. Then for any fixed $K$, because the diameter of $B^K$ is $\gamma/L$, if $x\in X$ and $\xi\in B^K$ satisfy $H(x,\xi)+\gamma\leq 0$, we have
	\begin{displaymath}
	H(x,\xi')\leq H(x,\xi)+L\|\xi-\xi'\|_\infty\leq H(x,\xi)+\gamma\leq 0\qquad\forall\xi'\in B^K.
	\end{displaymath}
	Therefore, by the independence of $\{\xi^j\}_{j \in [N]}$,
	\begin{displaymath}
	\begin{aligned}
	\pSamp\big(\hat{X}_{N,\gamma}&\nsubseteq X^{\text{Fea}}\big)\leq\pSamp\Big(\bigcup_{K\in\prod_{i \in [d]} [M_i]}
	\big\{\{\xi^\sset:\sset\in[N]^d\}\cap B^K=\emptyset \big\}\Big)\\
	=&\pSamp\Big(\bigcup_{i\in [d]}\bigcup_{k \in [M_i]}\big\{\{\xi_i^j:j\in[N]\}\cap[\bar{\xi}^k_i-\gamma/2L,\bar{\xi}^k_i+\gamma/2L]=\emptyset\big\}\Big)\\
	\leq&\sum_{i \in [d]}\sum_{k \in [M_i]}\pSamp\Big(\bigcap_{j \in [N]} \big\{\xi_i^j\notin\big[\bar{\xi}^k_i-\gamma/2L,\bar{\xi}^k_i+\gamma/2L\big]\big\}\Big)\\
	=&\sum_{i \in [d]}\sum_{k \in [M_i]} \prod_{j \in [N]} \pxi\big(\xi_i^j\notin\big[\bar{\xi}^k_i-\gamma/2L,\bar{\xi}^k_i+\gamma/2L\big]\big)\\
	\leq&\sum_{i \in [d]}M_i\big(1-\eta(\gamma)\big)^N\\
	\leq&(dDL/\gamma)\big(1-\eta(\gamma)\big)^N.
	\end{aligned}
	\end{displaymath}
	\qed
	For any $\beta\in(0,1),$ with sample size\begin{displaymath}
	N\geq\frac{\log(dDL/\gamma)+\log(1/\beta)}{\eta(\gamma)},
	\end{displaymath}
	we have $\pSamp\big(\hat{X}_{N,\gamma}\subseteq X^{\text{Fea}}\big)\geq 1-\beta$.
	
	In some cases, $\eta(\gamma)=\Omega(\gamma)$, for example, when the density function of $\xi_i$ is bounded away from 0 in $\Xi_i$. Therefore, in these cases, to obtain a completely reliable solution with confidence at least $1-\beta$, we need to generate $O\big(\big[\log(d/\gamma)+\log(1/\beta)\big]/\gamma\big)$ samples.

	\subsection{Two-stage stochastic linear programs with random right-hand side}\label{rhssec}
	
	Now we consider a two-stage stochastic linear programs with only random right-hand side, i.e., a problem of
	the form \eqref{real} - \eqref{recourse}, having $W(\xi)\equiv W,T(\xi)\equiv T$ independent of $\xi$. Recall that $H(x,\xi)$
	is defined in \eqref{hdef}. Note that if $h(\xi)$ is monotone in $\xi$, this case satisfies the monotonicity assumption of Example
	\ref{eg:monotone}, and so the theory in Section \ref{sec:prodmarg} applies. In this section we show that similar results
	can be obtained without imposing Assumption \ref{assum:prod} or the monotonicity assumption on $h(\xi)$.
	We find that applying a $\gamma$-``padding'' to the feasibility constraints of the original SAA
	formulation (\ref{stdsaa1}) is sufficient to obtain completely reliable solutions with high confidence.  
	
	Specifically, for $\gamma \in [0,\bar{\gamma}]$ define the SAA problem:
	\begin{equation}\label{modapprox2}
	\begin{aligned}
	\min \ &N^{-1}\sum_{j \in [N]} F(x,\xi^j)\\
	\text{s.t.} \ &H(x,\xi^j)+\gamma\leq 0,\quad j \in [N]\\
	&x\in X.
	\end{aligned}
	\end{equation}
	Let $\tilde{X}_{N,\gamma}$ denote the feasible region of the above problem. 
	\begin{thm}\label{prob2}
		The feasible region $\tilde{X}_{N,\gamma}$ of (\ref{modapprox2}) satisfies\begin{displaymath}
		\pSamp(\tilde{X}_{N,\gamma}\subseteq X^{\text{Fea}})\geq 1-n^W(1-\tilde{\eta}(\gamma))^N,
		\end{displaymath}
		where $n^W$ and $\tilde{\eta}(\gamma)$ are constants defined in (\ref{nwdef}) and (\ref{tetadef}), respectively. 
	\end{thm}
	\noindent {\it Proof}
	Let $E^W$ denote the set of all extreme points of $\bar{P}:=\{\alpha \in\bR_+^{m_2}:e_\cset^T\alpha=1,W^T\alpha=0\}$.
	By duality, $H(x,\xi)\leq 0$ is equivalent to $\alpha^T\big(h(\xi)-Tx\big) \leq 0$
	for all $\alpha \in E^W$, which implies that $H(x,\xi)\leq 0$ $\pxi$-almost surely if and only if
	\[ \alpha^TTx\geq\text{ess}\sup\big[\alpha^Th(\xi)\big] \]
	for all $\alpha \in E^W$. Similarly, $H(x,\xi^j)+\gamma\leq 0$ implies that\begin{displaymath}
	\alpha^TTx \geq \alpha^T\big[h(\xi^j)+\gamma e_\cset\big]=\alpha^Th(\xi^j)+\gamma
	\end{displaymath}
	for all $\alpha \in E^W$ and $j \in [N]$.
	Therefore, by independence of $\{\xi^j\}_{j \in [N]},$\begin{align*}
	\pSamp\big(\tilde{X}_{N,\gamma}\nsubseteq X^{\text{Fea}}\big)\leq&~\pSamp\Big(\bigcup_{\alpha \in E^W}\bigcap_{j \in [N]}\big\{\alpha^Th(\xi^j)+\gamma\leq\text{ess}\sup\big[\alpha^Th(\xi)\big] \big\}\Big)\\
	\leq&\sum_{\alpha \in E^W}\prod_{j \in [N]}\pxi\big(\big\{\alpha^Th(\xi^j)+\gamma\leq\text{ess}\sup\big[\alpha^Th(\xi)\big]\big\}\big)\\
	\leq &~n^W\big(1-\tilde{\eta}(\gamma)\big)^N
	\end{align*}
	where
	\begin{equation}\label{nwdef}
	n^W=|E^W|,
	\end{equation}
	\begin{displaymath}
	\tilde{\eta}(\alpha,\gamma):=\pxi\big(\alpha^Th(\xi)\geq\text{ess}\sup\big[\alpha^Th(\xi)\big]-\gamma\big) \quad (>0)
	\end{displaymath}
	and\begin{equation}\label{tetadef}
	\tilde{\eta}(\gamma):=\min_{\alpha \in E^W}\tilde{\eta}(\alpha,\gamma).
	\end{equation}
	\qed
	
	For any $\epsilon\in(0,1)$, since $n^W=|E^W|\leq\frac{m_2^{n_2+1}}{(n_2+1)!}\leq m_2^{n_2+1}$, with sample size\begin{displaymath}
	N\geq\frac{(n_2+1)\log(m_2)+\log(1/\epsilon)}{\tilde{\eta}(\gamma)},
	\end{displaymath}
	we have $\pSamp\big(\tilde{X}_{N,\gamma}\subseteq X^{\text{Fea}}\big)\geq 1-\epsilon$.
	
	Even though $\tilde{\eta}(\gamma)$ is the minimum value of potentially exponentially many $\tilde{\eta}_i(\gamma)$'s,
	$\tilde{\eta}(\gamma)=\Omega(\gamma)$ in some cases. For example, if $h(\xi)$ follows a uniform distribution
	$U([0,1]^{m_2})$ then $\tilde{\eta}(\gamma)=\gamma$ if $\gamma\in(0,1)$.

	
	\section{Numerical tests}\label{numerical}
	\subsection{A two-stage resource planning problem}
	We test the standard and padded SAA approaches for solving a two-stage resource planning (TRP) problem. This problem is inspired by a
	problem in \cite{luedtke2014branch}. The problem consists of a set of resources (e.g., server types), denoted by $i\in
	[n]$, which can be used to meet demands of a set of customer types, denoted by $k\in[m]$. The problem is stated as:
	\begin{equation}\label{testpr}
	\min_{x\in\bR^{n-p}_+ \times\bZ_+^{p}} c^Tx+\bE Q(x,\bm{\xi}),
	\end{equation}
	where for a fixed $\xi=(q,\rho,\mu,\lambda)$,
	\begin{displaymath}
	\begin{aligned}
	Q(x,\xi)=&\min_{y\in\bR^{n\times m}_+} &&q^Ty\\
	&\text{s.t.} &&\sum_{k \in [m]} y_{ik}\leq\rho_i x_i, && i \in [n],\\
	& && \sum_{i \in [n]} \mu_{ik}y_{ik}\geq\lambda_k, && k \in [m].
	\end{aligned}
	\end{displaymath}
	Here $c_i$ represents the unit cost of resource $i\in [n]$. For $i \in [n]$, variable $x_i$ represents the amount of resource $i$ to
	purchase and for $i \in [n], k \in [m]$ variable $y_{ik}$ represents the amount of resource $i$ allocated to customer type $k$ after observing
	the uncertainty $\xi$.  Parameters $q,\rho,\mu,\lambda$ are random vectors, where $q_{ik}$ represent the unit cost of
	allocating resource $i \in [n]$ to customer type $k
	\in [m]$, $\rho_i$ represents the utilization rate of resource $i
	\in [n]$, $\mu_{ik}$ represents the service rate of resource $i \in [n]$ for customer type $k \in [m]$ and $\lambda_k$
	represents the demand of customer type $k \in [m]$.
	
	\subsection{Test instances}\label{subsec:TestInst}
	We apply SAA approaches to the TRP problem on several instances where the first-stage
	variables are continuous $(p=0)$ or pure integer $(p=n)$. For the continuous first-stage case, we consider instances
	with $n\in\{10,20,40\}$, $m\in\{10,40\}$ and $N\in\{100,500,1000\}$. For the pure integer first-stage case, we consider
	instances with $n\in\{5,10,20\}$, $m\in\{10,40\}$ and $N\in\{50,100,500\}$. We generate $c$ following the scheme of
	\cite{luedtke2014branch}. To generate the scenario data $(q^j,\rho^j,\mu^j,\lambda^j)$, we generate the original
	distribution of $(q,\rho,\mu,\lambda)$ following \cite{liu2016decomposition}. We then truncate the original distribution
	such that the support of each entry of $(q,\rho,\mu,\lambda)$ is restricted to $[a-4\sigma,a+4\sigma]$, where $a$ and
	$\sigma$ are the mean and the standard deviation of that entry, respectively. We generate
	$(q^j,\rho^j,\mu^j,\bar{\lambda}^j)$ from the truncated distribution and then set $\lambda^j=\bar{\lambda}^j/10$. We use a smaller demand here in order to make the discrete version of the problem more distinct from the continuous version. Note that $q$ is a function of $\rho$ in \cite{liu2016decomposition}. Therefore, we define $\bm{\xi}$ as $(\rho,\mu,\lambda)$. In this case, $\bm{\xi}$ has a bounded support which is equal to the product of its marginal supports.
	\subsection{Recourse likelihood of SAA solutions}\label{subsec:reclike}
	\begin{table}
		\begin{center}
			\caption{Results for TRP problems with continuous first-stage variables}\label{contable}
			\begin{tabular}{ @{\extracolsep{\fill}} ccccr@{,}r}
				\toprule
				$(n,m,N)$ & \multicolumn{2}{c}{SAA objective value} & \multicolumn{3}{c}{$1-$Recourse likelihood (est.)}\\
				\midrule\midrule\addlinespace
				& Mean & \multicolumn{1}{c}{[Min,Max]} & Mean & \multicolumn{2}{c}{[Min,Max]}\\
				& (95\% C.I.) & \multicolumn{1}{c}{} & (95\% C.I.) & \multicolumn{2}{c}{}\\
				\cmidrule(lr){2-3}\cmidrule(lr){4-6}\addlinespace
				(10,10,\phantom{0}100)  & 	118.0$\pm 1.1$ & [113.8,124.3] & 2.38\%$\pm0.83\%$ & [0.10\% & 5.71\%]\\
				(10,10,\phantom{0}500)  & 	120.8$\pm 0.8$ & [118.2,124.3] & 0.57\%$\pm0.18\%$ & [0.05\% & 1.43\%]\\
				(10,10,1000) & 				121.9$\pm 0.6$ & [119.6,124.3] & 0.31\%$\pm0.09\%$ & [0.05\% & 0.74\%]\\
				(10,40,\phantom{0}100)  & 	434.3$\pm 4.9$ & [421.7,459.8] & 1.46\%$\pm0.51\%$ & [0.01\% & 3.55\%]\\
				(10,40,\phantom{0}500)  & 	442.3$\pm 2.9$ & [432.8,458.7] & 0.35\%$\pm0.13\%$ & [0.01\% & 1.01\%]\\
				(10,40,1000) & 				447.0$\pm 3.1$ & [439.1,462.5] & 0.17\%$\pm0.06\%$ & [0.00\% & 0.42\%]\\
				(20,10,\phantom{0}100)  & 	113.9$\pm 0.8$ & [109.6,116.5] & 3.55\%$\pm1.07\%$ & [1.37\% & 9.67\%]\\
				(20,10,\phantom{0}500)  & 	118.0$\pm 0.8$ & [115.4,121.8] & 0.58\%$\pm0.19\%$ & [0.08\% & 1.46\%]\\
				(20,10,1000) & 				118.9$\pm 0.6$ & [117.2,121.6] & 0.30\%$\pm0.09\%$ & [0.07\% & 0.68\%]\\
				(20,40,\phantom{0}100)  & 	454.3$\pm 3.0$ & [441.4,467.0] & 1.75\%$\pm0.57\%$ & [0.19\% & 3.90\%]\\
				(20,40,\phantom{0}500)  & 	460.1$\pm 1.9$ & [452.6,466.4] & 0.61\%$\pm0.18\%$ & [0.15\% & 1.57\%]\\
				(20,40,1000) & 				464.8$\pm 2.6$ & [455.8,476.6] & 0.30\%$\pm0.11\%$ & [0.02\% & 0.91\%]\\
				(40,10,\phantom{0}100)  & 	113.3$\pm 1.0$ & [109.1,117.1] & 3.41\%$\pm1.19\%$ & [0.53\% & 10.48\%]\\
				(40,10,\phantom{0}500)  & 	116.1$\pm 0.6$ & [114.0,118.5] & 0.72\%$\pm0.21\%$ & [0.14\% & 1.57\%]\\
				(40,10,1000) & 				117.3$\pm 0.4$ & [115.5,118.7] & 0.31\%$\pm0.07\%$ & [0.12\% & 0.77\%]\\
				(40,40,\phantom{0}100)  & 	461.2$\pm 2.9$ & [449.7,472.9] & 3.00\%$\pm0.74\%$ & [0.53\% & 6.09\%]\\
				(40,40,\phantom{0}500)  & 	473.1$\pm 1.8$ & [466.3,481.1] & 0.49\%$\pm0.12\%$ & [0.15\% & 0.98\%]\\
				(40,40,1000) & 				475.8$\pm 1.8$ & [470.2,486.9] & 0.28\%$\pm0.07\%$ & [0.03\% & 0.59\%]\\
				\bottomrule
			\end{tabular}
		\end{center}
	\end{table}
	\begin{table}
		\begin{center}
			\caption{Results for TRP problems with pure integer first-stage variables}\label{inttable}
			\begin{tabular}{ @{\extracolsep{\fill}} ccccr@{,}r}
				\toprule
				$(n,m,N)$ & \multicolumn{2}{c}{SAA objective value} & \multicolumn{3}{c}{$1-$Recourse likelihood (est.)}\\
				\midrule\midrule\addlinespace
				& Mean & \multicolumn{1}{c}{[Min,Max]} & Mean & \multicolumn{2}{c}{[Min,Max]}\\
				& (95\% C.I.) & \multicolumn{1}{c}{} & (95\% C.I.) & \multicolumn{2}{c}{}\\
				\cmidrule(lr){2-3}\cmidrule(lr){4-6}\addlinespace
				(\phantom{0}5,10,\phantom{0}50)  & 	137.9$\pm 1.1$ & [133.8,142.5] & 3.06\%$\pm1.14\%$ & [0.18\% & 9.85\%]\\
				(\phantom{0}5,10,100)  & 			140.3$\pm 1.2$ & [136.6,146.9] & 1.20\%$\pm0.40\%$ & [0.01\% & 3.20\%]\\
				(\phantom{0}5,10,500) & 			141.9$\pm 0.9$ & [139.3,147.4] & 0.42\%$\pm0.13\%$ & [0.01\% & 0.94\%]\\
				(\phantom{0}5,40,\phantom{0}50)  & 	470.2$\pm 6.7$ & [450.5,500.2] & 3.37\%$\pm1.67\%$ & [0.06\% & 10.67\%]\\
				(\phantom{0}5,40,100)  & 			473.2$\pm 4.8$ & [457.9,498.0] & 1.91\%$\pm0.79\%$ & [0.06\% & 4.92\%]\\
				(\phantom{0}5,40,500)  & 			487.3$\pm 3.8$ & [473.8,499.3] & 0.35\%$\pm0.15\%$ & [0.03\% & 1.04\%]\\
				(10,10,\phantom{0}50)  & 			117.0$\pm 1.4$ & [112.1,125.3] & 3.64\%$\pm1.38\%$ & [0.11\% & 10.38\%]\\
				(10,10,100)  & 						118.0$\pm 1.1$ & [113.9,124.4] & 2.39\%$\pm0.83\%$ & [0.10\% & 6.04\%]\\
				(10,10,500) & 						120.8$\pm 0.8$ & [118.3,124.3] & 0.54\%$\pm0.17\%$ & [0.06\% & 1.37\%]\\
				(10,40,\phantom{0}50)  & 			428.6$\pm 6.0$ & [406.4,464.0] & 3.31\%$\pm1.65\%$ & [0.01\% & 13.84\%]\\
				(10,40,100)  & 						434.3$\pm 4.9$ & [421.7,459.9] & 1.46\%$\pm0.51\%$ & [0.01\% & 3.56\%]\\
				(10,40,500) & 						442.3$\pm 2.9$ & [432.8,458.7] & 0.35\%$\pm0.13\%$ & [0.01\% & 1.00\%]\\
				(20,10,\phantom{0}50)  & 			112.8$\pm 1.2$ & [107.3,116.8] & 6.09\%$\pm2.43\%$ & [1.72\% & 18.71\%]\\
				(20,10,100)  & 						113.9$\pm 0.8$ & [109.6,116.5] & 3.46\%$\pm1.02\%$ & [1.40\% & 9.30\%]\\
				(20,10,500) & 						118.0$\pm 0.8$ & [115.4,121.9] & 0.57\%$\pm0.18\%$ & [0.09\% & 1.41\%]\\
				(20,40,\phantom{0}50)  & 			450.8$\pm 3.2$ & [440.4,463.3] & 2.66\%$\pm0.89\%$ & [0.59\% & 6.55\%]\\
				(20,40,100)  & 						454.3$\pm 3.0$ & [441.4,467.0] & 1.74\%$\pm0.57\%$ & [0.18\% & 3.94\%]\\
				(20,40,500) & 						460.1$\pm 1.9$ & [452.6,466.4] & 0.60\%$\pm0.18\%$ & [0.15\% & 1.57\%]\\
				\bottomrule
			\end{tabular}
		\end{center}
	\end{table}
	
	We summarize the standard SAA experiment results on the TRP instances with continuous first-stage variables and pure integer first-stage
	variables, respectively. Both experiments are conducted to observe how the SAA objective value and the recourse likelihood of the optimal solution is influenced by the problem size and sample size. For each combination of
	$(n,m,N)$, the same problem is solved using SAA with 20 different samples. Table \ref{contable} reports some statistics
	of the solutions generated for the TRP problem with continuous first-stage variables, and Table \ref{inttable} gives the
	same for the integer first-stage case. In particular, we report the SAA objective value and an estimate of the recourse likelihood $\phi(\hat{x}_N)$ of the solution $\hat{x}_N$, i.e., the
	probability the solution has a feasible recourse action. We actually display one minus the recourse likelihood (i.e.,
	the estimated probability that the solution does not have a feasible recourse action)
	as these numbers are easier to distinguish from each other. The recourse likelihood is estimated using an independent sample of 100,000
	iid scenarios. For each combination of $(n,m,N)$, we report in each row the range of the SAA objective value and
	recourse likelihood over the 20 samples, as well as an approximate 95\%-level confidence interval for the means of the SAA objective value and recourse likelihood.
	
	We first discuss the results for continuous first-stage variables. As a lower bound of the optimal value $v^*$, $\bE
	\hat{v}_N$ is monotonically increasing in the sample size $N$. In our experiments, we can observe that the sample means
	do increase as $N$ increases. The mean of violation probability is approximately proportional to $n/N$. Contrary to our
	theory, the violation probability does not show a clear connection to the number of second-stage variables in our
	experiment. Specifically, for fixed $n$ and $N$, the solutions obtained in instances with $m=40$ do not have significantly
	higher violation probability than those obtained with $m=10$, despite there being four times as many second-stage
	decision variables in these instances. This indicates that either these test instances have special structure not
	captured by our theory, or potentially that our analysis is not tight in terms of dependence on number of second-stage
	variables.
	
	We observe similar results in the experiments for the pure integer first-stage case. Compared with the continuous case,
	we observe a slightly lower violation probability for the same problem size and sample size. In most cases, a sample
	size $N=1000$ is still not enough for obtaining a potentially completely reliable solution (i.e., at least one of the
	100,000 scenarios used for assessing recourse feasibility was violated). 
	
	\subsection{Obtaining completely reliable solutions for the TRP problem via padded SAA}\label{subsec:testpaddedSAA_mon}
	We test the padded SAA approach for solving the TRP problem with continuous first-stage variables. In this problem,
	since the support $\bm{\xi}$ is equal to the product of its marginal supports, we apply a padding in feasibility
	constraints as described in Section \ref{hcsec}. Moreover, the TRP problem satisfies the monotonicity assumption in
	Example \ref{eg:monotone} with $\bm{\xi}=(-\rho,-\mu,\lambda)$. Therefore, we only need to pad on the most ``dominating" scenario. Specifically, the following problem is equivalent to the $\gamma$-padded SAA problem \eqref{modapprox}:
	
	\begin{equation}\label{TRP:padded}
	\begin{aligned}
	&\min_{(x,y,z)\geq 0} &&c^Tx+\frac{1}{N}\sum_{j \in [N]} (q^j)^Ty^j\\
	&\text{s.t.} &&\sum_{k \in [m]} y^j_{ik}\leq\rho^j_i x_i, && i \in [n],~j \in [N],\\
	& && \sum_{i \in [n]}\mu^j_{ik}y^j_{ik}\geq\lambda^j_k, && k \in [m],~j \in [N],\\
	&&&\sum_{k \in [m]} z_{ik}\leq\rho^{\min}_i x_i-\gamma, && i \in [n],\\
	&&&\sum_{i \in [n]} \mu^{\min}_{ik}z_{ik}\geq\lambda^{\max}_k+\gamma, && k \in [m].
	\end{aligned}
	\end{equation}
	Here $\rho^{\min}$ is defined as the componentwise minimum of $\{\rho^1,\ldots,\rho^N\}$, and $\mu^{\min}$ and
	$\lambda^{\max}$ are defined likewise. Note that the optimal objective value of \eqref{TRP:padded} is a piecewise linear convex function of $\gamma$. The padded problem \eqref{TRP:padded} is an LP with almost the same problem size as the original SAA problem.
	
	We summarize the padded SAA experiment results on the TRP instances with continuous first stage variables. We consider
	instances with $(n,m)$ fixed to be $(10,10)$ but with varying sample size $N\in\{100,500,1000\}$ and padding level
	$\gamma\in\{0,0.1,\ldots,1.9,2\}$. For each combination of $(N,\gamma)$, we solve the padded SAA problem
	\eqref{TRP:padded} with 20 different samples. For each $N$, we plot in Figures \ref{pSAA100}-\ref{pSAA1000} the average
	optimal objective value of \eqref{TRP:padded} and the fraction of solutions obtained from solving \eqref{TRP:padded}
	that are completely reliable over the 20 samples as $\gamma$ increases. To check if the solution is completely reliable, we check if the solution has a feasible recourse action for the ``hardest" scenario $\xi^{\sup}:=(\rho^{\inf},\mu^{\inf},\lambda^{\sup})$ where $\rho^{\inf}$ is the componentwise essential infimum of the random vector $\rho$, and $\mu^{\inf}$ as well as $\lambda^{\sup}$ is defined likewise.
	\begin{figure}[hbt!]
		\centering
		\includegraphics[width=0.6\linewidth]{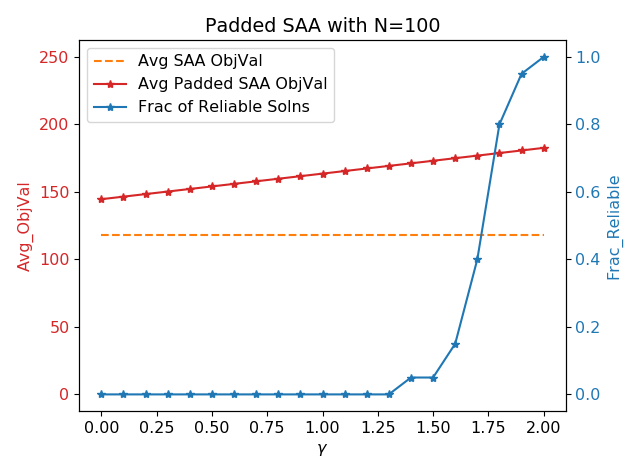}
		\caption{Average objective value and fraction of completely reliable solutions with $N=100$}
		\label{pSAA100}
	\end{figure}
	\begin{figure}[hbt!]
		\centering
		\includegraphics[width=0.6\linewidth]{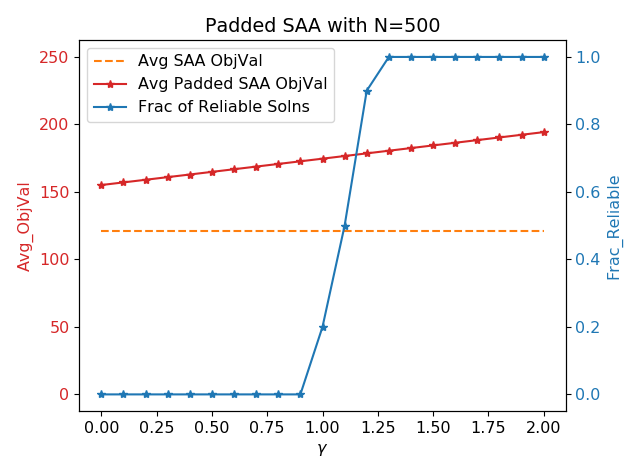}
		\caption{Average objective value and fraction of completely reliable solutions with $N=500$}
		\label{pSAA500}
	\end{figure}
	\begin{figure}[hbt!]
		\centering
		\includegraphics[width=0.6\linewidth]{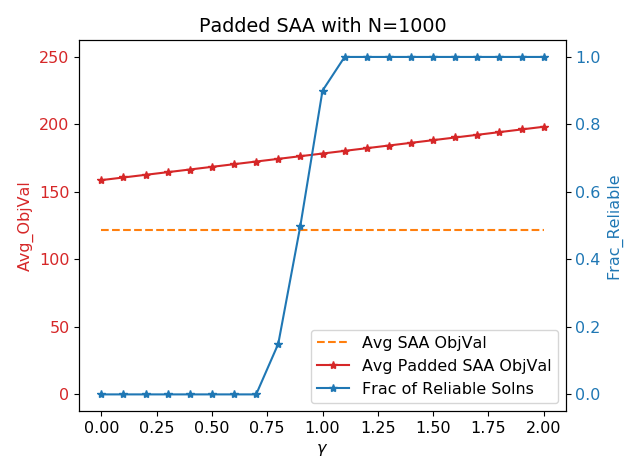}
		\caption{Average objective value and fraction of completely reliable solutions with $N=1000$}
		\label{pSAA1000}
	\end{figure}
	
	We observe that although the optimal objective value of \eqref{TRP:padded} is in general piecewise linear convex in $\gamma$,
	the average objective value is actually linear in $\gamma$ for the TRP instances when $\gamma\in[0,2]$. The objective value of
	the padded SAA problem is in general higher than the SAA objective value even when the padding level $\gamma$ is set to
	be 0 because a feasibility constraint for the mixed scenario $\xi^{\max}:=(\rho^{\min},\mu^{\min},\lambda^{\max})$ is
	added into the padded SAA problem. After reaching a certain threshold, e.g., $\gamma=0.7$ for $N=1000$, the fraction of
	completely reliable solutions obtained from different samples quickly increases from 0 to 1 as $\gamma$ increases.
	These results indicate that padded SAA indeed can be used to obtain completely reliable solutions.
	We also observe that in this example completely reliable
	solutions are significantly more costly than solutions obtained with standard SAA, which are not completely reliable but
	have high recourse likelihood. For example, with $N=1000$, the lowest cost completely reliable solution obtained is 44.6\% more costly than the average cost of the solution to the SAA with
	$N=1000$ scenarios, solutions which had an estimated recourse likelihood of 99.69\% on average (see Table \ref{contable}
	in Section \ref{subsec:reclike}). Thus, for this example we
	observe that if the decision-maker is willing to accept a small possibility that there is no feasible recourse,
	significant cost savings can be achieved in the (highly likely) case that a feasible recourse action does exist.
	
	\subsection{Test of padded SAA on a non-monotone variant of the TRP problem}\label{subsec:testpaddedSAA_nonmon}
	
	We next test the feasibility of the constraint-generation algorithm proposed in Section \ref{subsubsec:solving_padded}. We consider a non-monotone variant of the TRP problem, where the demand $\lambda_k$ of customer type $k$ is replaced by a linear factor model: $\lambda_k=\sum_{q\in[l]}a_{qk}\tau_q+h_k$ for $k\in[m]$,
	where $\{\tau_q\}_{q \in [l]}$ are random factors and
	$\{a_{qk}\}_{q \in [l]}$ and $h_k$ are deterministic parameters. We created test instances by generating the parameters
	$\{h_k\}_{k\in[m]}$ as independent realizations of a $N(11,6.25)$ random variable and the parameters
	$\{a_{qk}\}_{q\in[l],k\in[m]}$ as independent realizations of a $N(0,1)$ variable. Random variables
	$\{\tau_q\}_{q\in[l]}$ follow the uniform distribution $U([-1,1])$ and are mutually independent. The other random variables $(q,\rho,\mu)$ follow the same distribution as described in Section \ref{subsec:TestInst}, and $\bm{\xi}$ is defined as $(\rho,\mu,\tau)$.
	
	\begin{table}[htb]
		\begin{center}
			\caption{Strength of the two MILP formulations for solving \eqref{eq:sep}}\label{table:4MILP}
			\begin{tabular}{ @{\extracolsep{\fill}} ccccccccccccccccc}
				\toprule
				$(m,l)$ & \multicolumn{4}{c}{Proposition \ref{prop:separation}} & \multicolumn{4}{c}{ Proposition \ref{prop:separation+}}\\
				\cmidrule(lr){2-5}\cmidrule(lr){6-9}
				& $t_{\text{IP}}$(s) & gap & \# nodes & gap$_{\text{LP}}$ & $t_{\text{IP}}$(s) & gap & \# nodes & gap$_{\text{LP}}$\\
				\midrule
				(10,\phantom{0}5) & $>$600 & 32.7\% & $>1.18\times10^6$  & 73.1\% & 0.02 & 0\% & 0 & 0\%\\
				(10,10) & $>$600 & 26.6\% & $>4.12\times 10^5$ & 32.6\% & 0.08 & 0\% & 0 & 0\%\\
				(10,20) & $>$600 & 17.9\% & $>2.15\times 10^5$ & 20.4\% & 0.38 & 0\% & 0 & 0\%\\
				(40,20) & $>$600 & 33.1\% & $>7214$ & 33.6\% & 0.42 & 0\% & 0 & 0\%\\
				(40,40) & $>$600 & 15.3\% & $>2512$ & 15.4\% & 1.29 & 0\% & 0 & 0\%\\
				\bottomrule
			\end{tabular}
		\end{center}
	\end{table}
	
	We first test the ability of the MILP formulations from Propositions
	\ref{prop:separation} and \ref{prop:separation+}  to solve problem \eqref{eq:sep}. Observe that $H(x,\xi)$ is
	monotone in $(-\rho,-\mu)$, i.e., $H(x,\rho,\mu,\lambda)\leq H(x,\rho',\mu',\lambda)$ if $(\rho,\mu)\geq(\rho',\mu')$.
	Therefore, we can fix $\rho=\rho^{\min}$ and $\mu=\mu^{\min}$ in \eqref{eq:sep}. Note that the problem has fixed
	recourse with $\rho$ and $\mu$ fixed. 
	We test the MILP
	formulations on instances
	with $n=10$, $N=100$, and $(m,l)\in\{(10,5),(10,10),(10,20),(40,20),(40,40)\}$. For each
	combination of $(m,l)$, we first generate a candidate solution $\hat{x}$ by solving the SAA problem \eqref{saadef}, and
	then use the two different formulations to solve the associated constraint generating problem \eqref{eq:sep} with a 10
	minute time limit. 
	The results are presented in Table \ref{table:4MILP}, where `$t_{IP}$(s)' denotes the MILP solution time in seconds,
	`gap' denotes the optimality gap obtained at the time limit, `\# nodes' denotes the number of branch-and-cut nodes
	explored, and `gap$_{\text{LP}}$' denotes the LP relaxation gap. We observe that the MILP formulation from Proposition
	\ref{prop:separation+} has an extremely strong LP relaxation value and hence is solved for all test instances at the root node without branching. 
	On the other hand, the
	formulation from Proposition \ref{prop:separation} which does not exploit the fixed recourse structure explores a very
	large number of nodes and cannot be solved within the time limit for even a single instance at this size. We conclude
	from this experiment that when the montonicity assumption is not satisfied the algorithm proposed in Section
	\ref{subsubsec:solving_padded} has potential to be practical in the case of fixed recourse, but without fixed recourse
	would be limited to very small instances.

	We next apply the constraint-generation algorithm in Section \ref{subsubsec:solving_padded}  to solve the padded SAA
	problem \eqref{modapprox}, using the MILP formulation
	from  Proposition \ref{prop:separation+} to solve the subproblem \eqref{eq:sep}. We test on non-monotone TRP instances
	with $n\in\{10,40\}$, $m\in\{10,40\}$, $l=m$ and $N=1000$. For each combination of $(n,m)$, we report in Table
	\ref{table:nonmon} the total solution time, in seconds, used by the proposed algorithm (`Soln time(s)'), time spent on
	solving subproblems \eqref{eq:sep} (`MILP time(s)'), and the number of MILPs solved within the algorithm (`\# MILPs
	solved'), averaged over $11$ instances with $\gamma\in\{0,0.1,\ldots,0.9,1 \}$. We observe that only a very small number
	of ``mixed scenarios" $\xi^J$ are generated by the constraint generation algorithm. 
	These results indicate that the method is
	practical for instances with fixed recourse. 
	Note that  the vast majority of the solution time in these reported results is on solving the padded SAA master problem. While we do not pursue this, Benders
	decomposition could be applied to potentially accelerate this part of the method.
	
	\begin{table}
		\begin{center}
			\caption{Padded SAA results for non-monotone TRP Problems}\label{table:nonmon}
			\begin{tabular}{ @{\extracolsep{\fill}} cccc}
				\toprule
				$(n,m)$ & Soln time(s) & MILP time(s) & $\#$ MILPs solved\\
				\midrule\midrule\addlinespace
				(10,10) & \phantom{00}4 & \phantom{0}0.3 & 2\\
				(10,40) & \phantom{0}16 & \phantom{0}2.9 & 2\\
				(40,10) & \phantom{0}29 & \phantom{0}1.6 & 3\\
				(40,40) & 157 & 28.0 & 2\\
				\bottomrule
			\end{tabular}
		\end{center}
	\end{table}
	\section{Conclusion}
	
	We have presented some results on SAA for solving two-stage stochastic programs without relatively complete
	recourse. Our first results consider two-stage stochastic linear programs, and indicate that the probability the SAA
	solution has recourse likelihood less than $1-\epsilon$ convergences to zero exponentially fast. We obtained exact
	convergence results in terms of obtaining a completely reliable and near-optimal solution in the case that the feasible
	region is finite. Finally, we analyzed the use of ``padded'' SAA problems to obtain solutions that are completely
	reliable in cases when the feasible region is not finite. Numerical tests demonstrated empirically the relationship
	between the sample size and violation probability of the solutions obtained via SAA.
	
	As mentioned in Section \ref{saasec}, for two-stage stochastic linear programs, a sample size proportional to $n_1n_2$
	is sufficient for obtaining a solution with high recourse likelihood. We, however, observed in Section \ref{numerical}
	that the sample size required to obtain solutions with high recourse likelihood was independent of the number of second
	stage decision variables. Thus, it is an open question whether the theoretical results can be improved, or whether there
	exists problems for
	which the required sample size is $O(n_1n_2)$.
	
	\appendix
	\section*{Appendix}
	\subsection*{Proof of Corollary \ref{Nbound}}
	For any $\nu\in(0,1)$,
	\begingroup
	\allowdisplaybreaks
	\begin{align*} N\geq&\frac{1}{1-\nu}\Big[\frac{1}{\epsilon}\log\frac{1}{\beta}+\frac{n_1}{\epsilon}\bigl(\log \frac{1}{\nu\epsilon}+\log(\frac{m_1}{n_1}+2)+1\bigr)+\\
	&\frac{2n_1(n_2+1)}{\epsilon}\bigl(\log (\frac{m_2}{n_2+1})+1\bigr)+n_1\Big]\\
	\Rightarrow (1-\nu)N\geq&\frac{1}{\epsilon}\log\frac{1}{\beta}+\frac{n_1}{\epsilon}\bigl(\log \frac{1}{\nu\epsilon}+\log(\frac{2n_1+m_1}{n_1})+1\bigr)+\frac{2n_1(n_2+1)}{\epsilon}\bigl(\log (\frac{m_2}{n_2+1})+1\bigr)+n_1\\
	\Rightarrow N\geq&\frac{1}{\epsilon}\log\frac{1}{\beta}+\frac{n_1}{\epsilon}(\log \frac{1}{\nu\epsilon}+\log n_1-1+\frac{\nu N\epsilon}{n_1})+\frac{n_1}{\epsilon}\bigl(\log(\frac{2n_1+m_1}{n_1^2})+2\bigr)+\\
	&\frac{2n_1(n_2+1)}{\epsilon}\bigl(\log (\frac{m_2}{n_2+1})+1\bigr)+n_1\\
	\Rightarrow N\geq&\frac{1}{\epsilon}\log\frac{1}{\beta}+\frac{n_1}{\epsilon}\log N+\frac{n_1}{\epsilon}\bigl(\log(\frac{2n_1+m_1}{n_1^2})+2\bigr)+\\
	&\frac{2n_1(n_2+1)}{\epsilon}\bigl(\log (\frac{m_2}{n_2+1})+1\bigr)+n_1\\
	\Rightarrow\log\beta\geq& n_1\log N+n_1\bigl(\log(\frac{2n_1+m_1}{n_1^2})+2\bigr)+2n_1(n_2+1)\bigl(\log (\frac{m_2}{n_2+1})+1\bigr)+\\
	&\epsilon n_1-\epsilon N\\
	\Rightarrow\log\beta\geq& n_1\log N+n_1\log(2n_1+m_1)+2n_1(n_2+1)\log m_2+\\
	&\epsilon n_1-\epsilon N-2\log(n_1!)-2n_1\log\bigl((n_2+1)!\bigr)\\
	\Rightarrow \beta\geq&\frac{N^{n_1}(2n_1+m_1)^{n_1}m_2^{2n_1(n_2+1)}e^{-\epsilon(N-n_1)}}{(n_1!)^2\bigl((n_2+1)!\bigr)^{2n_1}}\\
	\geq& \binom{N}{n_1}\frac{1}{n_1!}(2n_1+m_1)^{n_1}\Bigl(\frac{m_2^{n_2+1}}{(n_2+1)!}\Bigr)^{2n_1}(1-\epsilon)^{N-n_1}.
	\end{align*}
	\endgroup
	The third inequality can be justified by observing that $\frac{\nu N\epsilon}{n_1}\geq 1+\log\frac{\nu N\epsilon}{n_1}$. The fifth inequality can be justified by observing that $\log (k!)\geq k(\log k -1)$ for $k\geq 1$. The result follows by setting $\nu=1/2$.
	
	\subsection*{Proof of Proposition \ref{lpprop}}
	We only need to prove that $H(x,\xi)$ is Lipshitz continuous in
	$(W(\xi),T(\xi),h(\xi))$ for all $x\in X$ under infinity (matrix) norm. Then the result
	follows from the assumption that $(W(\xi),T(\xi),h(\xi))$ is Lipshitz continuous in $\xi$.
	
	For fixed $x\in X$, assume $(y^*,\eta^*)$ and $(y',\eta')$ are optimal solutions of (\ref{hdef}) for
	$M:=(\randw,\randt,\randh):=(\randw(\xi),\randt(\xi),\randh(\xi))$ and
	$M':=(\randw',\randt',\randh'):=(\randw(\xi'),\randt(\xi'),\randh(\xi'))$, respectively. Consider the different conditions.\begin{enumerate}
		\item If the problem has only right-hand side randomness, then\begin{displaymath}
		\begin{aligned}
		\eta^*=&\min_y\{\max_{i \in \cset} \{\randh_i - \randw_iy-\randt_ix\} : Dy \geq d - Cx \}\\
		\leq&\max_{i \in \cset} \{\randh_i -\randw_iy'-\randt_ix\}\\ 
		\leq &\max_{i \in \cset} \{\randh'_i-\randw_iy'-\randt_ix\}+\|M-M'\|_\infty\|\left(
		\begin{array}{c}
		0\\
		0\\
		1
		\end{array}
		\right)\|_\infty\\
		\leq&\eta'+\|M-M'\|_\infty.
		\end{aligned}
		\end{displaymath}
		\item If the problem has fixed recourse and $X$ is bounded, then \begin{displaymath}
		\begin{aligned}
		\eta^*=&\min_y\{\max_{i \in \cset} \{\randh_i - \randw_iy-\randt_ix\}:Dy\geq d-Cx \}\\
		\leq&\max_{i \in \cset} \{\randh_i -\randw_iy'-\randt_ix\}\\
		\leq &\max_{i \in \cset} \{\randh'_i-\randw_iy'-\randt'_ix\}+\|M-M'\|_\infty\|\left(
		\begin{array}{c}
		0\\
		x\\
		1
		\end{array}
		\right)\|_\infty\\
		\leq&\eta'+R\|M-M'\|_\infty
		\end{aligned}
		\end{displaymath}
		for some $R>0$ since $X$ is bounded.
		\item If $\{(x,y):x\in X,Dy\geq d-Cx \}$ is bounded, then similarly\begin{displaymath}
		\begin{aligned}
		\eta^*=&\min_y\{\max_{i \in \cset} \{\randh_i - \randw_iy-\randt_ix\}:Dy\geq d-Cx \}\\
		\leq&\max_{i \in \cset} \{\randh_i -\randw_iy'-\randt_ix\}\\
		\leq &\max_{i \in \cset} \{\randh'_i-\randw'_iy'-\randt'_ix\}+\|M-M'\|_\infty\|\left(
		\begin{array}{c}
		y'\\
		x\\
		1
		\end{array}
		\right)\|_\infty\\
		\leq&\eta'+R\|M-M'\|_\infty.
		\end{aligned}
		\end{displaymath}
		for some $R>0$ since $\{(x,y):x\in X,Dy\geq d-Cx \}$ is bounded.
	\end{enumerate}
	
	So for each case, there exists $R>0$ such that
	\begin{displaymath}
	H(x,\xi)-H(x,\xi')=\eta^*-\eta'\leq R\|M-M'\|_\infty.
	\end{displaymath}
	Similarly, $H(x,\xi')-H(x,\xi)\leq R\|M-M'\|_\infty$.
	Therefore,
	\begin{displaymath}
	|H(x,\xi)-H(x,\xi')|\leq R\|M-M'\|_\infty,
	\end{displaymath}
	which implies that $H(x,\xi)$ is a Lipschitz continuous function in $(W(\xi),T(\xi),h(\xi) )$ for all $x\in X$ under infinity norm. Therefore, Assumption \ref{assum:lip} holds.
	
	\subsection*{Proof of Proposition \ref{prop:separation}}
	
	First observe that the linear program \eqref{hdef} is always feasible. Thus, 
	\begin{displaymath}
	\begin{aligned}
	H(x,\xi)
	=&\max_{\alpha} &&\alpha^T(\randh(\xi)-\randt(\xi)\hat{x})+\beta^T(d - C\hat{x}),\\
	&\text{s.t. }&&\randw(\xi)^T\alpha + D^T\beta =0,\\
	&&&e^T\alpha=1,~\alpha\geq 0,~\beta\geq 0,
	\end{aligned}
	\end{displaymath}
	where we adopt the convention that if the dual linear program is infeasible then the optimal value is defined to be
	$-\infty$.
	Thus,
	\begin{alignat*}{2}
	\max \{ H(\hat{x},\xi^\sset): \sset \in [N]^d \} =&
	\max_{\alpha,\beta, \sset}\ &&\alpha^T\big(\randh(\xi^\sset)-\randt(\xi^\sset)\hat{x}\big)+\beta^T(d-C\hat{x}) ,\\
	&\text{s.t.}&&\alpha^T\randw(\xi^\sset) + \beta^T D=0,\\
	&&&e^T\alpha=1,~\alpha\geq 0,~\beta \geq 0, \sset \in [N]^d 
	\end{alignat*}
	Next, introduce binary variables $\delta_{qj}$ for $q \in [d]$ and $j \in [N]$, where $\delta_{qj}=1$ implies that
	$\sset_q=j$. This leads to the  mixed-integer nonlinear program:
	\begin{equation}\label{sep:formulation1}
	\begin{aligned}
	&\max_{\alpha,\beta,\xi,\delta}\ &&\alpha^T(\randh(\xi)-\randt(\xi)\hat{x}) + \beta^T(d - C\hat{x})\\
	&\text{s.t. }&&\xi_q=\sum_{j \in [N]}\xi^j_q\delta_{qj},&&q \in [d],\\
	&&&\sum_{j \in [N]}\delta_{qj}=1,&&q \in [d],\\
	&&& \alpha^T \randw(\xi) + \beta^TD =0,\\
	&&&e^T\alpha=1,\\
	&&&\alpha\geq 0,~\beta \geq 0, ~\delta\in\{0,1\}^{d\times N}.
	\end{aligned}
	\end{equation}
	Using the assumptions that $\randw(\xi)$, $\randt(\xi)$ and $\randh(\xi)$ are linear in $\xi$, problem \eqref{sep:formulation1} can be written as the following mixed-integer bilinear program:
	\begin{alignat}{3}
	&\max_{\alpha,\beta,\xi,\delta}\ &&\alpha^T\Bigl(\overline{H}\xi-\sum_{k \in [n_1]}\randt^k\xi\hat{x}_k \Bigr) + \beta^T(d - C\hat{x}) \label{eq:sep2obj}\\
	&\text{s.t. }&&\xi_q=\sum_{j \in [N]}\xi^j_q\delta_{qj},&& q \in [d], \label{eq:sep2xi} \\
	&&&\sum_{j \in [N]}\delta_{qj}=1,&&q \in [d], \label{eq:sep2delt} \\
	&&&\alpha^T\randw^k\xi +\beta^T D^k =0,&&k \in [n_2], \label{eq:sep2dual} \\
	&&&e^T\alpha=1, \label{eq:sep2norm} \\
	&&&\alpha\geq 0,~\beta\geq 0,~\delta\in\{0,1\}^{d\times N}. \label{eq:sep2sign}
	\end{alignat}
	We next use \eqref{eq:sep2xi} to substitute out the variables $\xi_q$ in this formulation. Observe that one this is
	done, the only nonlinear terms are of the form $\alpha_p\delta_{qj}$ for $p \in \cset, q \in [d], j \in [N]$.
	Thus, introduce new variables $z_{pqj}$ to represent this product for each $p \in \cset$, $q \in
	[d]$, and $j \in [N]$. 
	Using constraint \eqref{eq:sep2norm} we derive the linear constraints: 
	\begin{equation}
	\label{eq:zlin1}
	\sum_{p \in \cset} z_{pqj} \Bigl(= \sum_{p \in \cset} \alpha_p\delta_{qj} \Bigr) = \delta_{qj}, \quad q \in
	[d],~j \in [N].
	\end{equation}
	Using constraints \eqref{eq:sep2delt} we derive the linear constraints: 
	\begin{equation}
	\label{eq:zlin2}
	\sum_{j \in [N]} z_{pqj} \Bigl(= \sum_{j \in [N]} \alpha_p\delta_{qj} \Bigr) = \alpha_p, \quad p \in \cset, \ q \in [d].
	\end{equation}
	Observe that constraints \eqref{eq:zlin1},\eqref{eq:zlin2} together with \eqref{eq:sep2sign} and \eqref{eq:sep2delt} are
	sufficient to imply $z_{pqj} = \alpha_p \delta_{qj}$. Indeed for any fixed $q$, suppose $j^*_q$ is the index such that
	$\delta_{qj^*_q}=1$. Then \eqref{eq:zlin1} implies that $z_{pqj}=0=\alpha_p\delta_{qj}$ for any $j \neq j^*_q$ and all $p$,
	and \eqref{eq:zlin2} implies that $\sum_{j \in [N]} z_{pqj} = z_{pqj^*_q} = \alpha_p = \alpha_p \delta_{qj^*_q}$ also for
	all $p$. Also note that constraints \eqref{eq:sep2norm}, \eqref{eq:zlin1} and \eqref{eq:zlin2} imply\begin{displaymath}
	\sum_{j\in[N]}\delta_{qj}=\sum_{j\in[N]}\sum_{p\in I}z_{pqj}=\sum_{p\in I}\Big(\sum_{j\in[N]}z_{pqj}\Big)=e^T\alpha=1
	\end{displaymath}
	for $q\in[d]$. Therefore, constraints \eqref{eq:sep2delt} are redundant when constraints \eqref{eq:sep2norm}, \eqref{eq:zlin1} and \eqref{eq:zlin2} are present.
	Thus the mixed-integer bilinear program \eqref{eq:sep2obj}-\eqref{eq:sep2sign} can be reformulated as the MILP
	given in Proposition \ref{prop:separation}  by using \eqref{eq:sep2xi} to substitute out the variables $\xi$, using $z_{pqj}$ to replace the bilinear terms
	$\alpha_p\delta_{qj}$, and adding the constraints \eqref{eq:zlin1}-\eqref{eq:zlin2}.
	
	\subsection*{Proof of Proposition \ref{prop:separation+}}
	Note that $H(\hat{x},\cdot)$ is convex under the assumptions and the maximum of a convex function can only be attained at extreme points. Therefore, we can reformulate \eqref{eq:sep} as\begin{equation}\label{eq:sep+}
	\max\Big\{H(\hat{x},\xi):\xi\in\prod_{q\in[d]}\{\xi^{\min}_{q},\xi^{\max}_{q} \} \Big\}.
	\end{equation}
	Introduce binary variables $\delta_{q1}$ and $\delta_{q2}$ for $q\in[d]$, where $\delta_{q1}=1$ implies that $\xi_q=\xi^{\min}_{q}$ and $\delta_{q2}=1$ implies that $\xi_q=\xi^{\max}_{q}$. We can then rewrite \eqref{eq:sep+} as a mixed-integer bilinear program:
	\begin{equation}\label{sep+:formulation1}
	\begin{aligned}
	&\max_{\alpha,\beta,\xi,\delta}\ &&\alpha^T\Bigl(\overline{H}\xi-\sum_{k \in [n_1]}\randt^k\xi\hat{x}_k \Bigr) + \beta^T(d - C\hat{x})\\
	&\text{s.t. }&&\xi_q=\xi^{\min}_q\delta_{q1}+\xi^{\max}_q\delta_{q2},&&q \in [d],\\
	&&&\delta_{q1}+\delta_{q2}=1,&&q \in [d],\\
	&&& \alpha^T \randw + \beta^TD =0,\\
	&&&e^T\alpha=1,\\
	&&&\alpha\geq 0,~\beta \geq 0, ~\delta\in\{0,1\}^{d\times 2}.
	\end{aligned}
	\end{equation}
	Similar to Proposition \ref{prop:separation}, we introduce variables $z_{pqj}$ to represent $\alpha_p\delta_{qj}$ for $p\in I,q\in[d],j\in\{1,2\}$. Then problem \eqref{eq:sep+} can be written as the following mixed-integer linear program:
	\begin{align}
	&\max_{\alpha,\beta,\delta,z}\ &&\sum_{p\in I}\sum_{q\in[d]}\Big[\Big(\overline{H}_{pq}-\sum_{k\in[n_1]}\hat{x}_k\randt^k_{pq}\Big)(\xi_q^{\min}z_{pq1}&&\hspace{-3mm}+\xi_q^{\max}z_{pq2}) \Big]+ \beta^T(d - C\hat{x})\label{sep+:formulation2}\\
	&\text{s.t. }&&\sum_{p\in I}z_{pqj}=\delta_{qj},&&q\in[d],~j\in\{1,2\},\label{sep+:sum_z}\\
	&&&z_{pq1}+z_{pq2}=\alpha_p,&&p\in I,~q\in[d],\label{sep+:alpha}\\
	&&& \alpha^T \randw + \beta^TD =0,\label{sep+:dual}\\
	&&&e^T\alpha=1,\label{sep+:norm}\\
	&&&\alpha\geq 0,~\beta \geq 0,~z\geq 0, ~\delta\in\{0,1\}^{d\times 2}.\label{sep+:sign}
	\end{align}
	Finally, we apply the reformulation-linearization technique to strengthen the MILP formulation
	\eqref{sep+:formulation2}-\eqref{sep+:sign}. We introduce new variables $w_{pqj}\geq 0$ to represent the product
	$\beta_p\delta_{qj}$ for $p\in [m_2]\setminus I,q\in[d],j\in\{1,2\}$. Note that $w_{pqj}=\beta_p\delta_{qj}$ and $\alpha^T\randw+\beta^TD=0$ imply\begin{multline}\label{con:rlt}
	\sum_{p\in I}\randw_{pk}z_{pqj}+\sum_{p\in [m_2]\setminus I}D_{pk}w_{pqj}\Big(=\delta_{qj}\big(\sum_{p\in I}\randw_{pk}\alpha_p+\sum_{p\in [m_2]\setminus I}D_{pk}\beta_p \big)\Big)=0,\\
	k\in[n_1],q\in[d],j\in\{1,2\},
	\end{multline}
	and constraints $\delta_{q1}+\delta_{q2}=1$ for $q\in[d]$ imply\begin{equation}\label{con:beta_sum}
	w_{pq1}+w_{pq2}\Big(=\beta_p(\delta_{q1}+\delta_{q2})\Big)=\beta_p,\quad p\in[m_2]\setminus I,~q\in[d].
	\end{equation}
	Note that constraints \eqref{sep+:alpha}, \eqref{con:rlt} and \eqref{con:beta_sum} imply \eqref{sep+:dual}.
	Therefore, adding new variables $w\geq 0$ together with constraints \eqref{con:rlt} and \eqref{con:beta_sum} into the original MILP formulation yields a strengthened formulation in the lifted space, as the MILP given in Proposition \ref{prop:separation+}.

\bibliographystyle{ieeetr}
\bibliography{reference}   

%
%

\end{document}